\documentclass[11pt]{amsart}
\usepackage{enumerate}
\usepackage{amsopn,amsfonts,amsmath,amssymb}
\usepackage{tikz-cd}
\usepackage{array}	
\usepackage{color}
\newcommand{\indep}{\perp \!\!\! \perp}
\usepackage[colorlinks=true,linkcolor=red,urlcolor=blue,citecolor=gray]{hyperref}
\usepackage[top=1.20in, bottom=1.50in, left=1.25in, right=1.25in]{geometry}

\def\E{{\mathbb E} }	
\def\indic{{\rm {\large 1}\hspace{-2.3pt}{\large l}}}

\def\R{{\mathbb R}}
\def\Z{{\mathbb Z}}
\def\N{{\mathbb N}}

\def\C{{\mathbb C}}
\def\supp{{\rm supp}}

\def\Aff{{\rm Aff}}

\newcommand{\df}[2]{\frac{#1 }{#2}}

\newcommand{\abs}[1]{ \left|  #1\right| }
\newcommand{\mt}[1]{ \boldsymbol{ #1 } }

\newtheorem{theorem}{Theorem}

\newtheorem{proposition}{Proposition}

\newtheorem{corollary}{Corollary}
\newtheorem{assumption}{Restriction}
\newtheorem*{assumption*}{Restriction}
\newtheorem{definition}{Definition}
\newtheorem{example}{Example}

\newtheorem{remark}{Remark}

\usepackage[foot]{amsaddr}

\begin{document}
\title[Nonparametric classes for identification in random coefficients models]{Nonparametric classes for identification in random coefficients models when regressors have limited variation}

\author[Gaillac]{Christophe Gaillac$^{{(1)}}$}
\address{$ ^{(1)}$ Toulouse School of Economics, Universit\'e Toulouse Capitole, 1 esplanade de l'universit\'e, 31000 Toulouse, France.}
\email{\href{mailto:christophe.gaillac@tse-fr.eu}
{christophe.gaillac@tse-fr.eu}}

\author[Gautier]{Eric Gautier$^{(1)}$}
\email{\href{mailto:eric.gautier@tse-fr.eu}{eric.gautier@tse-fr.eu}}
\date{This version: \today.} 

\thanks{The authors acknowledge financial support from the grants ERC POEMH 337665 and ANR-17-EURE-0010. Christophe Gaillac thanks CREST/ENSAE where this research was partly conducted.}

\begin{abstract}
This paper studies point identification of the distribution of the coefficients in some random coefficients models with exogenous regressors when their support is a proper subset, possibly discrete but countable. We exhibit trade-offs between restrictions on the distribution of the random coefficients and the support of the regressors. We consider linear models including those with nonlinear transforms of a baseline regressor, with an infinite number of regressors and deconvolution, the binary choice model, and panel data models such as single-index panel data models and an extension of the Kotlarski lemma.
\end{abstract}

\maketitle

\section{Introduction}\label{s1} 
Random coefficients models are used to incorporate multiple sources of unobserved heterogeneity in modelling various economic behaviors. In these models, the researcher can be interested in more than the average of the vector of random coefficients and even might want to recover its distribution (see, \emph{e.g.}, \cite{berry2009nonparametric,fox2017note,GK,hoderlein2014triangular,masten2017random}). Imposing parametric assumptions on the law of the random coefficients is a widely used approach in the analysis of random coefficients models (see, \emph{e.g.},  \cite{berry1995automobile,Train}) but can seriously drive the results (see  \cite{breunig2019,HS}). Economic theory rarely motivates such restrictions. When the coefficients are
random with a law in a nonparametric class and independent from the regressors, identification often requires the regressors to have a support which is the whole space. 
However, in applications regressors may only have limited variation. 

The main contribution of this paper is to provide identification of the distribution of the coefficients in some random coefficients models with regressors independent from the random coefficients, namely exogenous ones, when their support is a proper subset, possibly discrete but countable. Estimation results along those lines for the linear model are available in \cite{estimation}. We provide a general identification strategy which is then used to study linear models  including those with nonlinear transforms of a baseline regressor, with an infinite number of regressors  and deconvolution, the binary choice model, and panel data models such as single-index panel data models and an extension of the Kotlarski lemma. We show in these models that the support of the regressors can be discrete if we maintain integrability assumptions on the distribution of the coefficients. Using several examples, we illustrate the trade-off between the variation the regressors and the restrictions on the nonparametric class of distributions of random coefficients we consider. 

The paper is organized as follows. Section \ref{sec:not} introduces the notations and the nonparametric identification strategy. In Section \ref{sec:Fourier}, we provide the main identification results for the various linear random coefficients models. Section \ref{s:BC} provides results on the random coefficients binary model and Section \ref{sec:panel_all} on some panel data models with random coefficients. The appendix provides complements on the tools used for identification and 
 the proofs.

\section{Preliminaries}\label{sec:not} 
Bold letters are used for vectors, capital letters for indeterminates of polynomials or random variables/vectors. For a real number $r$, $\boldsymbol{r}_c$ is the vector, the dimension of which will be clear from the text, where each entry is $r$. 
The notations $\N$ and $\N_0$ are used for the positive and nonnegative integers, for $p\in\N$ $[p]=\{1,\dots,p\}$, $\R_+=\{x\in\R:\ x>0\}$, $\mathbb{C}(X)$ denotes  the set of rational functions, and $\indic\left\{\cdot\right\}$ the indicator function. $\mathbb{K}[Z_1,\dots,Z_p]$ is the ring of polynomials of $p$ variables with coefficients in the ring $\mathbb{K}$. Let $S \subseteq \R^p$, $C(S)$ and $C^{\infty}(S)$ be the continuous and infinitely differentiable functions at every point in $S$ with values in $\mathbb{C}$.
For $\mathcal{S}\subseteq\C^{p}$, $\mathcal{A}(\mathcal{S})$ and 
$\mathcal{H}^{\infty}(\mathcal{S})$ are the analytic and bounded analytic functions on $\mathcal{S}$. Quasi-analytic classes of functions on $\R^{p}$ are  vector spaces of complex valued functions in $C^{\infty}(\R^{p})$ characterized by the knowledge of their derivatives at $\bold{0}_c$. Unlike real analytic functions, the Taylor series of such functions do not need to converge. 
 The notation $|\cdot|_q$ for $q\in[1,\infty]$ stands for the $\ell_q$ norm of a vector with components in $\C$ while $L^{q}(S,\mu)$ for $q\in[1,\infty]$ are the $q$ integrable functions on $S$ with respect to the measure $\mu$ and the norms $\|\cdot\|_{L^{q}(S,\mu)}$. When $\mu$ is the Lebesgue measure we drop it from the notation. 
$\mathbb{S}^p$ 
is the unit sphere in $\mathbb{R}^{p+1}$. 
$(\boldsymbol{e}_j)_{j=1}^{p}$ is the canonical basis of $\R^{p}$. 
For $\boldsymbol{\beta}\in\C^{p}$, 
$k\in \N_0$, and $\boldsymbol{m}\in\mathbb{N}_0^{p}$, denote by $\boldsymbol{m}!=\prod_{j=1}^{p}\boldsymbol{m}_j!$, $|\boldsymbol{m}|=\sum_{j=1}^{p}\boldsymbol{m}_j$, $\boldsymbol{\beta}^{\boldsymbol{m}}=\prod_{j=1}^{p}\boldsymbol{\beta}_j^{\boldsymbol{m}_j}$, and $|\boldsymbol{\beta}|^{\boldsymbol{m}}=\prod_{j=1}^{p}|\boldsymbol{\beta}_j|^{\boldsymbol{m}_j}$. 
$GL(\R^{p})$ are the invertible $p\times p$ matrices with real coefficients. 

 $\mathfrak{M}_c(\Omega)$, $\mathfrak{M}(\Omega)$, and $\mathfrak{M}_1(\Omega)$ are the sets of complex, nonnegative, and probability measures on a Borel measurable set $\Omega$. When $\Omega\subset\R^p$ is closed, these are measures over $\R^p$ with support in $\Omega$. For $\mu\in\mathfrak{M}_c(\R^{p})$ and $\boldsymbol{m}\in\N_0^{p}$, $|\mu|$ is its total variation, $s_{\mu}(\boldsymbol{m})=\int_{\R^{p}}\boldsymbol{x}^{\boldsymbol{m}}d\mu(\boldsymbol{x})$ the moments, and  by $s_{|\cdot|,|\mu|}(\boldsymbol{m})=\int_{\R^{p}}|\boldsymbol{x}^{\boldsymbol{m}}|d|\mu|(\boldsymbol{x})$ the absolute moments. Define
 \begin{align*}
\mathfrak{M}_c^{*}\left(\R^{p}\right)&=\left\{\mu\in\mathfrak{M}_c\left(\R^{p}\right):\ \forall m \in \N, \  \sum_{j=1}^p s_{|\cdot|,|\mu|}(2m\boldsymbol{e}_j)<\infty\right\},\\ 
\mathfrak{M}^{*}\left(\R^{p}\right)&=\mathfrak{M}_c^{*}\left(\R^{p}\right)\cap\mathfrak{M}\left(\R^{p}\right),\ \mathfrak{M}_1^{*}\left(\R^{p}\right)=\mathfrak{M}^{*}\left(\R^{p}\right)\cap\mathfrak{M}_1\left(\R^{p}\right).
\end{align*} 
Let 
$\mathcal{M}$ be $\mathfrak{M}_c^*\left(\Omega\right)$ or $\mathfrak{M}^*\left(\Omega\right)$. A measure $\mu$ is said to be determinate in $\mathcal{M}$ if $\mu \in \mathcal{M}$ and if $\mu$ is uniquely determined in $\mathcal{M}$ by $\{s_{\mu}(\boldsymbol{m})\}_{\boldsymbol{m}\in\N_0^{p}}$. 

The Fourier transform of $\mu\in\mathfrak{M}_c(\R^{d})$ (resp. $f$ in $L^q(\mathbb{R}^{d})$ for $q=1,2$) is $\mathcal{F}\left[\mu\right]: \ \boldsymbol{x} \mapsto \int_{\mathbb{R}^{d}}e^{i\boldsymbol{b}^{\top}\boldsymbol{x}}d\mu(\boldsymbol{b})$ (resp. $\mathcal{F}\left[f\right]$). For a random vector $\boldsymbol{X}$, $\mathbb{P}_{\boldsymbol{X}}$ is its law, $F_{\boldsymbol{X}}$ its CDF,  $f_{\boldsymbol{X}}$ its density, 
$\varphi_{\boldsymbol{X}}=\mathcal{F}\left[\mathbb{P}_{\boldsymbol{X}}\right]$, and $\mathbb{S}_{\boldsymbol{X}}$ its support. 
$ \E_{\mathbb{P}}\left[ \cdot \right] $ is the expectation under $\mathbb{P}$, $\otimes$ the product of measures, and $\mathbb{P}_{Y|\boldsymbol{X}}(\cdot|\boldsymbol{x})$ for $\boldsymbol{x}\in\mathbb{S}_{\boldsymbol{X}}$ the conditional distribution of $Y$ given $\boldsymbol{X}=\boldsymbol{x}$.

This paper considers models of the form  
\begin{equation}\label{ev} 
\boldsymbol{Y} = \boldsymbol{v}(\boldsymbol{X},\boldsymbol{\Gamma}),\ \boldsymbol{X}\indep\boldsymbol{\Gamma},
\end{equation}
where $\boldsymbol{v}$ is a known measurable vector valued function, 
$\boldsymbol{Y}$ 
is the vector of outcomes, $\boldsymbol{\Gamma}\in\Gamma$ and $\boldsymbol{X}\in\mathcal{X}$ 
are vectors of unobserved and observed factors, and $\Gamma$ and $\mathcal{X}$ are the Euclidian space or the sphere. 
Everything in this paper holds if we impose independence given $\boldsymbol{Z}$, where $\boldsymbol{Z}$ is a random vector  from which we could have observations simultaneously with those of $\boldsymbol{Y}$ and $\boldsymbol{X}$ or which could be identifiable from a model for $\boldsymbol{X}$ obtained by a control function approach. The independence means that $\boldsymbol{X}$ is exogenous, it yields
$$\forall \boldsymbol{x}\in\mathbb{S}_{\boldsymbol{X}},\ \mathbb{P}_{\boldsymbol{Y}|\boldsymbol{X}}(\cdot|\boldsymbol{x})=\mathbb{P}_{\boldsymbol{v}(\boldsymbol{x},\boldsymbol{\Gamma})}.$$ 

In order to study identification, we proceed as if 
$(\mathbb{P}_{\boldsymbol{Y}|\boldsymbol{X}=\boldsymbol{x}})_{\boldsymbol{x}\in\mathbb{S}_{\boldsymbol{X}}}$ were known and denote the true law by $\mathbb{P}_{\boldsymbol{\Gamma}}^*$. We  maintain restrictions $\mathcal{R}$ on the  primitives $\mathbb{P}_{\boldsymbol{X}}\otimes \mathbb{P}_{\boldsymbol{\Gamma}}$ which  involve two conditions:  
(i) $V\subseteq \mathbb{S}_{\boldsymbol{X}}$ and (ii) $\mathbb{P}_{\boldsymbol{\Gamma}}\in\mathcal{P}$, for well chosen sets $V$ and $\mathcal{P}$.  We use the notation $(\mathbb{S}_{\mt{X}},\mathbb{P}_{\boldsymbol{\Gamma}})\in \mathcal{R}$.
Random coefficients models involve inner products of subvectors of $\boldsymbol{\Gamma}$ and $\boldsymbol{X}$, hence indices with random coefficients, and possibly additional idiosyncratic errors. The linear random coefficients model can be written as
\begin{equation}\label{eq:linear}
Y = \alpha+ \boldsymbol{\beta}^{\top}\boldsymbol{X},\ \boldsymbol{\Gamma}=(\alpha,\mt{\beta}^{\top})^{\top}\indep\boldsymbol{X}.
\end{equation} 
It is a baseline component of more complex random coefficient models. For example, \cite{Gautier2} shows how identification of a discrete choice model with random coefficients and a special regressor is amenable to identification of \eqref{eq:linear} while without special regressor requires a treatment like in Section \ref{s:BC}. Obtaining identification of $\mathbb{P}_{\boldsymbol{\Gamma}}^*$ in these non separable models with multiple unobservables without restricting $\mathcal{P}$ usually requires  that $\mathbb{S}_{\boldsymbol{X}}=\mathcal{X}$. For model \eqref{eq:linear} 
this can be shown by the Cramer-Wold Theorem. 
$\mathbb{S}_{\boldsymbol{X}}=\mathcal{X}$ is too demanding for a dataset. 
However, identification can hold when $\mathbb{S}_{\boldsymbol{X}}$ is a proper subset 
but $\mathcal{P}$ is a restricted nonparametric class. Because random coefficient models involve indices, we can study identification when $\mathbb{S}_{\boldsymbol{X}}$ is replaced by a convenient invertible affine transformation. So when the results in this paper depend on the order of the regressors, the order is irrelevant.  
Indeed, for all $\underline{\boldsymbol{x}}\in\R^p$ and $M\in GL(\R^p)$,
\begin{equation}\label{eq:repar}
\alpha+\boldsymbol{\beta}^{\top}\boldsymbol{X}=\alpha+\boldsymbol{\beta}^{\top}\underline{\boldsymbol{x}}+\boldsymbol{\beta}^{\top}M^{-1}M\left(\boldsymbol{X}-\underline{\boldsymbol{x}}\right)
\end{equation}
and there is a one to one relation between $\mathbb{P}_{\alpha+\boldsymbol{\beta}^{\top}\underline{\boldsymbol{x}},\left(M^{-1}\right)^{\top}\boldsymbol{\beta}}$ and $\mathbb{P}_{\alpha,\boldsymbol{\beta}}$.

Let  $ \mathcal{T}^{\boldsymbol{v}}_{\mathbb{S}_{\boldsymbol{X}}}$ map $\mathbb{P}_{\boldsymbol{\Gamma}}$ to the collection  of direct image measures 
$\left(\mathbb{P}_{\boldsymbol{v}(\boldsymbol{x},\boldsymbol{\Gamma})}\right)_{\mt{x}\in\mathbb{S}_{\boldsymbol{X}}}$.  
$\mathbb{P}_{\boldsymbol{\Gamma}}^*$ is identified under $\mathcal{R}$ if $\mathcal{T}^{\boldsymbol{v}}_{\mathbb{S}_{\boldsymbol{X}}}$ is injective.
The proofs that we give often rely on  the diagram 
\begin{equation}\label{ediag}
\begin{tikzcd}[column sep=3.5pc]
\mathbb{P}_{\boldsymbol{\Gamma}}\in \mathcal{P} \quad \arrow{r}{\mathcal{T}_{\mathbb{S}_{\boldsymbol{X}}}^{\boldsymbol{v}}} \quad\arrow{d}{\mathcal{G}\text{ (injective)}}& \left(\mathbb{P}_{\boldsymbol{v}(\boldsymbol{x},\boldsymbol{\Gamma})}\right)_{\mt{x}\in \mathbb{S}_{\boldsymbol{X}}} \arrow{d}{\Pi}  \\
\mathcal{G}\left[\mathbb{P}_{\boldsymbol{\Gamma}}\right]\in
\mathfrak{F}(S) 
\quad \arrow{r}{\text{restriction to }U} &
\quad \mathcal{G}\left[\mathbb{P}_{\boldsymbol{\Gamma}}\right]_{|U} 
\end{tikzcd}
\end{equation}
The choice of $\Pi$ depends on $\boldsymbol{v}$. 
By restricting $\mathcal{P}$, $\mathcal{G}[\mathbb{P}_{\boldsymbol{\Gamma}}]$ belongs to a class of functions $\mathfrak{F}(S)$ which also contains all differences between functions in $\mathcal{G}[\mathbb{P}_{\boldsymbol{\Gamma}}]$ (\emph{e.g.}, $\mathfrak{F}(S)$ could be a vector space), $\mathcal{G}$ is injective, and  
$U$ is a set of uniqueness of $\mathfrak{F}(S)$.
\begin{definition}
$U\subseteq S$ is a set of uniqueness for a vector space of functions $\mathfrak{F}(S)$ on $S$ if every function of $\mathfrak{F}(S)$ which is zero on $U$ is identically zero on $S$.
\end{definition}

\section{The linear random coefficients model}\label{sec:Fourier}
The model is \eqref{eq:linear}. It is a natural specification to account for heterogenous effects. It can be viewed as more general than the quantile regression. Indeed, when the conditional quantiles are strictly increasing, the quantile regression defines the same data generating process as a linear random coefficients model where the coefficients are functions of a scalar uniform distribution. 
The unobserved scalar uniform variable is a ranking variable. \eqref{eq:linear} allows for the coefficients to be a function of a multidimensional vector, possibly infinite. \cite{Gautier2} allows for extensive generalizations of \eqref{eq:linear} including ones involving nonlinear transforms of a baseline regressor. This paper considers different specifications. In Section \ref{s:nonlinear}, \eqref{eq:linear} is an approximation of the nonlinear model 
\begin{equation}\label{eq:NL}
\boxed{\quad Y=g(\boldsymbol{X},\boldsymbol{\Theta}),\ \boldsymbol{\Theta}\indep \boldsymbol{X},}
\end{equation}
where $g$ is unknown, $\boldsymbol{\Theta}$ has arbitrary, possibly infinite, dimension. Indeed, when $\mathbb{S}_{\boldsymbol{X}}$ is compact and almost surely in $\boldsymbol{\theta}\in\mathbb{S}_{\boldsymbol{\Theta}}$, $\boldsymbol{x} \mapsto g(\boldsymbol{x},\boldsymbol{\theta})\in L^2(\mathbb{S}_{\boldsymbol{X}})$, 
$(f_{j})_{j\in\Z}$ is a Riesz basis and $f_0=1$,  then 
\begin{equation}\label{finitesum}
g(\boldsymbol{x},\boldsymbol{\theta})= \gamma_{0}(\boldsymbol{\theta})+\sum_{j\in\Z\setminus \{0\}}\gamma_{j}(\boldsymbol{\theta})f_{j}(\boldsymbol{x}).
\end{equation}
Recall that a Riesz basis of a separable Hilbert space is the image of an orthonormal basis by a bounded invertible operator and that the coefficients $\gamma_{j}(\boldsymbol{\theta})$ are the inner products of $\boldsymbol{x}\mapsto g(\boldsymbol{x},\boldsymbol{\theta})$ with a uniquely defined biorthogonal system (see sections 1.7 and 1.8 in \cite{Young}). An approximation of \eqref{eq:NL}
is thus 
\begin{equation}\label{eq:SNL}
Y=\gamma_{0}(\boldsymbol{\Theta})+\sum_{j\in J}\gamma_{j}(\boldsymbol{\Theta})f_{j}(\boldsymbol{X}),\ \boldsymbol{\Theta}\indep \boldsymbol{X},
\end{equation}
where $J\subset\Z\setminus \{0\}$ is finite. 
Equation \eqref{eq:SNL} can be used for extrapolation of the conditional distribution of $Y$ given $\boldsymbol{X}=\boldsymbol{x}$ when the functions $f_{j}$ are analytic (see also \cite{Note}). Hence, it is possible to extrapolate not only conditional expectation functions but also any conditional quantiles even under a nonparametric specification with multiple unobservables entering in a non-additively separable way. 
In \eqref{eq:NL} or its approximation \eqref{eq:SNL}, the researcher can also be interested in the law of the random coefficients to eventually obtain the law of the elasticity $\partial_x g(x,\boldsymbol{\Theta})x/g(x,\boldsymbol{\Theta})$ or of the marginal effects $\partial_x g(x,\boldsymbol{\Theta})$ assuming they exist and $x$ belongs to the interior of $\mathbb{S}_{X}$.  
\begin{remark} A variation of \eqref{eq:SNL} consists in taking $J=\Z\setminus\{0\}$ but there exists $j_0\in \N$ such that $\boldsymbol{\gamma}_j(\boldsymbol{\Theta})$ for $|j|>j_0$ are deterministic and $\sum_{j\in\Z} \mathbb{E}\left[\boldsymbol{\gamma}_j(\boldsymbol{\Theta})^2\right] < \infty$. 
$\E\left[\gamma_j(\boldsymbol{\theta}) \right]$ are identified because $ \E\left[\gamma_j(\boldsymbol{\theta}) \right] = \E\left[ Y g_j(\boldsymbol{X})/f_{\boldsymbol{X}}(\boldsymbol{X}) \right]$, where 
$(g_j)_{j\in \Z}$ is the biorthogonal system associated to $(f_j)_{j\in \Z}$.  So this model is amenable to \eqref{eq:SNL} with an a priori unknown set $J$. 
 However, by taking $j_0$ to be the smallest $k$ such that 
 \begin{equation*} 
j_0=\text{argmin}_{k\in\N}\inf_{\mathbb{P}_{\boldsymbol{\Gamma}}}\int_{\R}\mathbb{E}\left[\left(\mathbb{P}\left(Y\le y\left|\boldsymbol{X}\right.\right)- \mathbb{P}\left(\sum_{|j|\leq k}\boldsymbol{\Gamma}_{j}f_j(\boldsymbol{X})\leq y\right)\right)^2\right]dy.
\end{equation*}
\end{remark}
Finally, we also consider in this section infinite dimensional linear random coefficients models. 
There, the random coefficients 
belongs to $\ell_2(\N_0)$ a.s. It is denoted by $\{\Gamma_m\}$ and is such that $\Gamma_0=\alpha$ and, for all $m\in\N$, $\Gamma_m=\beta_m$. The regressors $\{X_m\}$ belong to $\ell_2(\N)$ a.s. and the model can be written
\begin{equation}\label{eq:inf}
Y = \alpha + \sum_{m=1}^{\infty}\beta_m X_m,\  \{\Gamma_m\} \indep \{X_m\}. 
\end{equation} 
It can be related to the random coefficients functional linear regression (see \cite{wu2010varying})  
$$ Y = \int_{-1}^1 \pi(s) X(s) ds + U,$$ 
where a.s. $X\in L^2(-1,1)$,  $ \pi = \sum_{m=0}^{\infty} \Pi_{m} \nu_m\in L^2(-1,1)$ and   
$X_m : = \int_{-1}^1\nu_k(s)X(s) ds$ belong respectively to $\ell_2(\N_0)$ and $\ell_2(\N)$, and $(\nu_m)_{m\in\N_0}\in L^2(-1,1)^{\N_0}$.
Similar expressions are used for model  \eqref{eq:inf}. 

In this section, we consider different identifying restrictions which we denote by $\mathcal{R}_{L,j}$, for $j\in[7]$ and are defined below.
The main theorem of this section is the following. 
\begin{theorem}\label{tL}
$\mathbb{P}_{\boldsymbol{\Gamma}}^*$ in \eqref{eq:linear} is identified under either of $\mathcal{R}_{L,j}$ for $j\in\{1,2,3,6,7\}$. 
$\mathbb{P}_{\{\Gamma_m\}}^*$ in \eqref{eq:inf} is identified under either of $\mathcal{R}_{L,j}$ for $j\in\{4,5\}$.  Moreover, $\mathbb{P}_{\boldsymbol{\Gamma}}^*$ in \eqref{eq:linear} is not identified if $\mathcal{R}_{L,3}$ \eqref{RL2ii} holds but not $\mathcal{R}_{L,3}$ \eqref{RL2i}.
\end{theorem}

The proofs rely on diagram \eqref{ediag} with 
$\mathcal{G} = \mathcal{F}$ (sometimes extended to an open set of $\C^{p+1}$), $U = \{ (t,t\mt{x}), \ (t,\mt{x}) \in\R\times \mathbb{S}_{\boldsymbol{X}}\}$, and 
$$ \Pi: \ \left(\mathbb{P}_{(1,\boldsymbol{x}^{\top})\boldsymbol{\Gamma}}\right)_{\mt{x}\in\mathbb{S}_{\boldsymbol{X}}} \mapsto \left((t,\mt{x}) \in\R\times \mathbb{S}_{\mt{X}} \mapsto 
 \int_{\R^{p+1}} e^{it(1,\boldsymbol{x}^{\top})\boldsymbol{\gamma}}d\mathbb{P}_{\boldsymbol{\Gamma}}(\boldsymbol{\gamma})\right),$$ 
with the relevant modifications when dealing with sequences and model \eqref{eq:inf}. We denote by $\mathcal{P}(\Omega)$ (resp. $\mathcal{P}_c(\Omega)$) the set of (resp. complex) measures with support included in $\Omega$ which are determinate in $\mathfrak{M}^*(\Omega)$ (resp. $\mathfrak{M}_c^*(\Omega)$). Appendix \ref{sec:CD} reviews criteria for determinacy and classical examples are given in Remark \ref{r:CM}. The vector space spanned by $\mathcal{F}[\mathcal{P}(\Omega)]$ (resp. $\mathcal{F}[\mathcal{P}_c(\Omega)]$) is a quasi-analytic class of functions on $\R^{p}$. It is known that sets with a nonempty interior and, when $p=1$, sets which contain a bounded sequence of distinct points are sets of uniqueness. For analytic classes this is the Weierstrass theorem (see, \emph{e.g.}, Theorem 15.11 in  \cite{Rudin2}) and for quasi-analytic classes this follows by the arguments in the proof of Lemma 4.8 in \cite{Belisle} for the quasi-analytic class under consideration. We now give details on the restrictions. 

\subsection{Regressors with finite support and independence of the marginals of $\mathbb{P}_{\boldsymbol{\Gamma}}$,  deconvolution}\label{sec:indep}
This is a challenging situation in terms of limited variation of the regressors. We present restrictions in the spirit of \cite{Beran1}.  Let $\Omega_0\subseteq\mathbb{R}$ and, for all $k\in[p]$, $\Omega_k\subseteq\mathbb{R}$ be given closed sets. These sets account for possible prior information on the support of the marginals of $\mathbb{P}_{\boldsymbol{\Gamma}}$. When such an information is not available, we take the set to be $\R$. We consider either of the two following restrictions. 
\begin{assumption}	
\begin{enumerate}[$\mathcal{R}_{L,1}$ \textup{(}i\textup{)}]
		\item\label{RL1i}
		$\left\{\boldsymbol{0}_c,\boldsymbol{x}(1),\dots,\boldsymbol{x}(p)\right\} \subseteq\mathbb{S}_{\boldsymbol{X}}$ and 
		$(\boldsymbol{x}(1),\dots,\boldsymbol{x}(p))$ is an upper triangular matrix with nonzero diagonal elements;
		\item\label{RL1ii} $\mathcal{P}=\{\mathbb{P}_{\boldsymbol{\Gamma}}=\mathbb{P}_{\boldsymbol{\alpha}}\otimes\bigotimes_{k=1}^p\mathbb{P}_{\boldsymbol{\beta}_k}:\ \forall k\in[p],\ \mathbb{P}_{\boldsymbol{\beta}_k}\in\mathcal{P}(\Omega_k)\}$. 
	\end{enumerate}
\end{assumption}
\begin{assumption}	
\begin{enumerate}[$\mathcal{R}_{L,2}$ \textup{(}i\textup{)}]
		\item\label{RL1pi} $\mathcal{R}_{L,1}$ \eqref{RL1i} holds ;
		\item\label{RL1pii} $\mathcal{P}=\{\mathbb{P}_{\boldsymbol{\Gamma}}=\mathbb{P}_{\boldsymbol{\alpha}}\otimes\bigotimes_{k=1}^p\mathbb{P}_{\boldsymbol{\beta}_k}:\ \mathbb{P}_{\boldsymbol{\alpha}}\in\mathcal{P}(\Omega_0),\ \forall k\in[p-1],\ \mathbb{P}_{\boldsymbol{\beta}_k}\in\mathcal{P}(\Omega_k)\}$. 
	\end{enumerate}
\end{assumption}

The problem of deconvolution with two samples, one of the error and one of the sum of the signal and the error, fits model \eqref{eq:linear} with $p=1$ and $\mathbb{S}_{X}=\{0,1\}$ is a particular case (see \cite{Gautier2}).   

$\mathcal{R}_{L,2}$ \eqref{RL1pii} is in the spirit of the assumption used in \cite{deconv}. 
Under $\mathcal{R}_{L,1}$ \eqref{RL1ii}, $\varphi_{\alpha}$ can have zeros on an open set at the expense of a stronger assumption on $\mathbb{P}_{\beta}$. 
There are classical examples of characteristic functions with compact support (\emph{e.g.}, $t\mapsto (1-|t|^{r})\indic\{|t|\le 1\}$ for $0<r\le 1$). Remark 1 in \cite{kotlarski} relies on analyticity. \cite{Meister07} considers the estimation of compactly supported densities of $\beta$ which implies $\mathcal{R}_{L,1}$ \eqref{RL1ii} and $\varphi_{\beta}$ is analytic. 

Now on, we do not assume mutual independence of the random coefficients. This is important if $\boldsymbol{\Gamma}$ is of the form $\boldsymbol{\Gamma}(\boldsymbol{\Theta})$, where $\boldsymbol{\Theta}$ is a deep heterogeneity parameter $\boldsymbol{\Theta}$ as in  \eqref{eq:SNL}.

\subsection{$\mathbb{S}_{\boldsymbol{X}}$ is not in the zeros of a nonzero polynomial and nonlinear model}\label{s:nonlinear}
Let $\Omega\subseteq \mathbb{R}^{p+1}$ be a closed set and the researcher knows that $\mathbb{S}_{\boldsymbol{\Gamma}}\subseteq C$. 
\begin{assumption}
	\begin{enumerate}[$\mathcal{R}_{L,3}$ \textup{(}i\textup{)}]
		\item\label{RL2i}
		$\mathbb{S}_{\mt{X}}$ is not a subset of the zeros of a nonzero element of $\R[Z_1,\dots,Z_{p}]$;
		\item\label{RL2ii} $\mathcal{P}$ is the set of measures which are determinate in $\mathfrak{M}^*(\Omega)$.  
	\end{enumerate}
\end{assumption}
Lemma 2 in \cite{masten2017random} considers the case where $\Omega=\R^{p+1}$ and assumes $\mathbb{S}_{\boldsymbol{X}}$ contains an open ball which implies $\mathcal{R}_{L,3}$ \eqref{RL2i} (see also \cite{BM}). $\mathcal{R}_{L,3}$ \eqref{RL2ii} is much weaker than assuming that $\mathbb{S}_{\boldsymbol{\Gamma}}$ is compact as in \cite{BM}. By the usual properties of the Fourier transform, the vector space spanned by $\mathcal{F}[\mathcal{P}]$ is a quasi-analytic class. Because $\mathcal{A}(\R^{p+1})$ is a small subset of $C^{\infty}(\R^{p+1})$ (see Appendix \ref{sec:meager}), it is important in the analysis of the paper to use larger classes and allow for Fourier transforms to not be analytic. 
The next examples satisfy $\mathcal{R}_{L,3}$ \eqref{RL2i} (\emph{i.e.} are sets of uniqueness of $\mathbb{R}_d[Z_1,\dots,Z_{p}]$).
\begin{example}\label{P1}
\begin{enumerate}[\textup{(}i\textup{)}]
	\item $p=1$ and $|\mathbb{S}_{\boldsymbol{X}}|=\infty$ ;
	\item $p\ge2$ and $U_{p}\subseteq \mathbb{S}_{\boldsymbol{X}}$, where $U_{p}$ is defined recursively via
	$U_{1}$, such that $|U_{1}|=\infty$ and, for all $j=2,\dots,p$, 
	$U_{j}=\bigcup_{\boldsymbol{u}\in{U_{j-1}}}\{ (\boldsymbol{u}^{\top},v)^{\top},\ v\in {V}_{j}(\boldsymbol{u})\}$, where $|{V}_{j}(\boldsymbol{u})|=\infty$.
\end{enumerate}
\end{example}
This is a set of uniqueness because $P\in\mathbb{R}_d[Z_1,\dots,Z_{p}]$ can be written as $P(Z_1,\dots,Z_p)=\sum_{k=0}^dQ_k(Z_1,\dots,Z_{p-1})Z_p^k$, where $Q_k\in\mathbb{R}_d[Z_1,\dots,Z_{p-1}]$.  Examples of sets $U_{p}$ are $\prod_{k=1}^{p}V_k \subseteq \mathbb{S}_{\mt{X}}$, where, for all $k\in[p]$, $|V_k|=\infty$, the infinite fan $\{\boldsymbol{x}\in\R^2:\ \exists n\in\N, \ \boldsymbol{x}_2=n\boldsymbol{x}_1\}$,
and infinite staircase $\{\boldsymbol{x}\in\R^2:\ \boldsymbol{x}_2=\left \lceil{\boldsymbol{x}_1}\right \rceil\}$ (see \cite{BCR}). 

A consequence of the second statement in Theorem \ref{tL} is that $\mathbb{P}_{\boldsymbol{\Gamma}}^*$ in \eqref{eq:linear} when $\boldsymbol{X}=(X\ X^2)^{\top}$ under $\mathcal{R}_{L,3}$ \eqref{RL2ii} holds, even if $\mathbb{S}_X=\R$. Indeed, $\mathbb{S}_{\boldsymbol{X}}$ is included in the set of zeros of $Q(Z_1,Z_2) = Z_1^2 - Z_2$. Example \ref{P2} shows that we can handle other nonlinear transformations of a baseline variable which can have  discrete support. It applies to model \eqref{eq:SNL}.  

\begin{example}\label{P2}$\{(f_1(\boldsymbol{u}),\dots,f_p(\boldsymbol{u}))^{\top},\ \boldsymbol{u}\in U_q \}\subseteq \mathbb{S}_{\boldsymbol{X}}$, where $q\in\N$, $\mathcal{B}(\mathcal{S})$ is a vector space of functions on $\mathcal{S}\subseteq\C^q$ which is an algebra with respect to multiplication, the set of functions $\{f_j, j\in[p]\}$ of $\mathcal{B}(\mathcal{S})$ is such that if $P(f_1,\dots,f_p)=0$, with $P\in \R[Z_1,\dots,Z_p]$, then $P=0$, and $U_q$ a set of uniqueness of $\mathcal{B}(\mathcal{S})$ (see Section \ref{sec:tools} for classes of analytic and quasi-analytic functions).
\end{example}
Because $\mathcal{B}(\mathcal{S})$ is an algebra with respect to multiplication, $P(f_1,\dots,f_p)\in\mathcal{B}(\mathcal{S})$. Because $U_q$ is a set of uniqueness of $\mathcal{B}(\mathcal{S})$, if $P(f_1,\dots,f_p)(\boldsymbol{x})=0$ for all $\boldsymbol{x}\in U_q$ then $P(f_1,\dots,f_p)=0$, so $P=0$. Hence, $\mathcal{R}_{L,3}$ \eqref{RL2i} holds. Here, $p$ can be large and $U_q$ discrete. 
The last condition in Example \ref{P2} means that $\{f_j, j\in[p]\}$ is algebraically independent over $\R$. Clearly, algebraic independence over $k$ implies algebraic independence over $\R$ if $\R$ is a subfield of $k$.  We now give useful examples to apply Example \ref{P2}.

If $f_1(\boldsymbol{z})=\boldsymbol{z}_1,\dots,f_{p-1}(\boldsymbol{z})=\boldsymbol{z}_q$, $f_p$ is a real or complex analytic function of $\boldsymbol{z}_1,\dots,\boldsymbol{z}_{p-1}$  and  $\{f_j, j\in[p]\}$ is a set of algebraically independent over $\C$, $f_p$ is said to be transcendental in $\C$. If $f$ is a transcendental function in $\C$, then so is $1/f$ and its inverse on a domain where it is invertible. If $f$ is meromorphic in $\C$ then it is a transcendental function in $\C$ iff it is not a rational function (\emph{e.g.}, $\exp(z)/z$, $\sin(z)/(z-1)^2$, $\cosh(z)$, $\Gamma(z)$, $\zeta(z)$). Classical examples of transcendental functions include $z^{\pi}$, $\exp(z)$,  $z^z$, $z^{1/z}$, and $\log(z)$.

If $q\le p-1$ and $f_1(\boldsymbol{x})=\boldsymbol{x}_1,\dots,f_q(\boldsymbol{x})=\boldsymbol{x}_q$, then $\{f_j, j\in[p]\}$ is an algebraically independent set over $\C$ iff $\{f_j, j\in\{q+1,\dots,p\}\}$
 is algebraically independent over $\C(Z_1,\dots,Z_q)$. 
For example, $\{f_j, j\in[p-1]\}$ are algebraically independent over $\C(Z)$ if: 
\begin{enumerate}[\textup{(}i\textup{)}]
\item\label{casei} $f_j(z)=e^{\lambda_jz}$, where $(\lambda_j)_{j=1}^p$ are such that, if $(b_j)_{j\in J} \in \mathbb{N}_0^{p}$ is such that $ \sum_{j=1}^p b_j \lambda_j = 0$, then, for all $j=1,\dots,p$, $b_j=0$ (see proof of Proposition \ref{prop:basis}). 
\item\label{caseii} $f_j=e^{\varphi_j}$ and $\varphi_j\in \mathcal{A}_0(\mathcal{S})$, where $\mathcal{A}_0(\mathcal{S})\subseteq\mathcal{A}(\mathcal{S})$ and $\mathcal{S}\subseteq\C$ is a simply connected open subset which does not contain 0 but a set of uniqueness of $\mathcal{A}_0(\mathcal{S})$, such that, for $j=2,\dots,p$, $\lim_{x\to\infty}|\varphi_j(x)/\varphi_{j-1}(x)|=\infty$ and $\lim_{n\to\infty}|\varphi_1(x)/\log(x)|=\infty$. Thus, $\left\{\Psi,e^{\varphi_1\circ\Psi},\dots,e^{\varphi_p\circ\Psi}\right\}$ is algebraically independent over $\C$ on a connected open subset $\mathcal{S}_1\subseteq\C$ if $\Psi:\ \mathcal{S}_1\to\mathcal{S}$ is analytic.
\end{enumerate}
A way to achieve the condition in \eqref{casei} above is to rely on a transcendental number $r$ over $\Z$ (\emph{e.g.}, $e$ or  $\pi$). It means that, for all $P\in \Z[Z]$, $P(r)=0$ implies that $P=0$. The following result shows that we can work with the specification \eqref{eq:SNL}
with functions from a Riesz basis which are algebraically independent. \begin{proposition}\label{prop:basis}
	Let $T\in\R_+$, $\lambda_j = j +1/(5r^{|j|})$ and $f_j(z)=e^{i \pi \lambda_j z/T }$, for all $j\in\Z\setminus\{0\}$ and $r\in (1,\infty)$ is transcendental over $\Z$, and $f_0=1$.
	We have:
	\begin{enumerate}[\textup{(}P{6}.1\textup{)}]
		\item\label{P11} 
		$\{f_j:\ j\in \Z\setminus\{0\}\}$ is an algebraically independent family of functions over $\C(X)$; 		\item\label{P12} $\left(f_j\right)_{j\in \Z}$ is a Riesz basis of $L^2(-T,T)$.
	\end{enumerate}
\end{proposition} 
(P6.\ref{P11}) implies $\{f_j:\ j\in \Z\setminus\{0\}\}$ is an algebraically independent family of functions over $\R$ which is the condition in Example \ref{P2}. Here $q=1$ and $U_q$ should be such that $U_q\subseteq [-T,T]$. 

Consider model \eqref{eq:inf}.  $(\{e_{j,m}\})_{j\in\N}$ is now the canonical basis of $\ell_2(\N)$. Let $G$ be a vector subspace of $\ell_2(\N)$ and $\mu$ a probability measure on $\ell_2(\N)$, then we denote by $\Pi_{G*}\mu$ the projection of $\mu$ onto $G$.

\begin{assumption} There exists $n_0\in\N$ and $\{\underline{x}_m\}\in \ell_2(\N)$ such that, for all $n\ge n_0$, 
	\begin{enumerate}[$\mathcal{R}_{L,4}$ \textup{(}i\textup{)}]  		
	\item\label{RL3i}
	$ \{ (x_1,\dots,x_n) : \{x_m\}\in \mathbb{S}_{\{X_m\}},\ \forall  j> n, x_j = \underline{x}_j\}  $
		is not a subset of the zeros of a nonzero element of $\R[Z_1,\dots,Z_{n}]$; 
		\item\label{RL3ii}  $\mathcal{P}$ is the set of measures whose projections onto
			$G_{n+1} := \text{Span}(\{e_{1,m}\},\dots, \{e_{n+1,m}\})$ 
		are determinate in $\mathfrak{M}^*(\R^{n+1})$.
	\end{enumerate}
\end{assumption}
$\mathcal{R}_{L,4}$ \eqref{RL3i} allows discrete support of the regressors. 
Like in Section \ref{sec:indep}, less restrictive conditions on $\mathbb{S}_{\{X_m\}}$ can be obtained assuming independence between marginals of $\mathbb{P}_{\{\Gamma_m\}}$. An alternative restriction is: 
\begin{assumption}  
\begin{enumerate}[$\mathcal{R}_{L,5}$ \textup{(}i\textup{)}]  
		\item\label{RL4i} There exists a Gaussian measure $\mu$ on $\ell_2(\N)$ such that  
	$
	 \{\{u_m\}\in\R^{\N_0}:\ \forall n\in\N,\ u_n=u_0x_n,\ \{x_m\}\in \mathbb{S}_{\{X_m\}}\}  $ has positive $\mu$-measure; 
		\item\label{RL4ii}  $\mathcal{P}$ is the set of measures such that, for all $n\in\N$ and $\{f_m\}\in\ell_2(\N)\setminus\{\{0_m\}\}$, the  projections onto
		$G_{n+1} := \text{Span}(\{e_{1,m}\},\dots, \{e_{n,m}\}, \{f_m\})$ 
		are determinate in $\mathfrak{M}^*(\R^{n+1})$.
	\end{enumerate}
\end{assumption}
 
\subsection{Without restrictions on $\alpha$} 
$\mathcal{R}_{L,3}$ \eqref{RL2ii} places restrictions on $\mathbb{P}_{\alpha}$ which we entirely remove. This is important to model income or wealth and allow $\alpha$ to have, for example, a Pareto distribution which does not even belong to $\mathfrak{M}^*(\R)$.  
We can replace the role of $\alpha$ by $\boldsymbol{\beta}_k$ for $k\in[p]$ or a combination of the coefficients if, starting from \eqref{eq:repar},  we form, for example, 
\begin{equation*}
\frac{Y}{\boldsymbol{X}_k-\underline{\boldsymbol{x}}_k}=\boldsymbol{\beta}_k+\frac{\alpha}{\boldsymbol{X}_k-\underline{\boldsymbol{x}}_k}+\sum_{j\ne k}\boldsymbol{\beta}_j\frac{\boldsymbol{X}_j-\underline{\boldsymbol{x}}_j}{\boldsymbol{X}_k-\underline{\boldsymbol{x}}_k}
\end{equation*}
and there exists $\underline{\boldsymbol{x}}_k\notin \mathbb{S}_{\boldsymbol{X}_k}$. If the support of a subvector of dimension $k$ of $\boldsymbol{X}$ is $\R^k$ then we can similarly handle situations where there are no restrictions on the joint distribution of the corresponding coefficients and $\alpha$.  In this section, $\Omega\subseteq\R^p$ is a closed set 
 and it can be known that $\mathbb{S}_{\boldsymbol{\beta}}\subseteq\Omega$. 
 
 A key idea is to rely on the partial Fourier transform. For example, $\mathcal{F}[\mathbb{P}_{\boldsymbol{\Gamma}}]$ can be analyzed as the collection of functions $\mathcal{F}[\mathbb{P}_{\boldsymbol{\Gamma}}](t,\cdot)$ for $t\in\R$ and 
$$\forall t\in\R,\ \mathcal{F}\left[\mathbb{P}_{\boldsymbol{\Gamma}}\right](t,\star)=\mathcal{F}\left[\mathbb{P}_{\boldsymbol{\beta},t}\right](\star),\ \text{where}\ 
\mathbb{P}_{\boldsymbol{\beta},t}=\int_{\R}e^{ita}d\mathbb{P}_{\boldsymbol{\Gamma}}(a,\cdot).$$
Because the vector space spanned by $\mathcal{F}[\mathcal{P}_c(\Omega)]$ is a quasi-analytic class of functions on $\R^{p}$ and by the continuity of the Fourier transform at $0$, $\mathbb{P}_{\boldsymbol{\Gamma}}^*$ in \eqref{eq:linear} is identified under $\mathcal{R}_{L,6}$ below.  
\begin{assumption}\label{R5} \begin{enumerate}[$\mathcal{R}_{L,6}$\ \textup{(}i\textup{)}]
	\item\label{it:C1}  $\mathbb{S}_{\mt{X}}$ has a nonempty interior or, when $p=1$,  $\mathbb{S}_{X}$ contains a bounded sequence of distinct points;
	 \item\label{it:C2}  $\mathcal{P}$ restricts $\mathbb{P}_{\boldsymbol{\Gamma}}$ so that  $\mathbb{P}_{\boldsymbol{\beta},t}\in\mathcal{P}_c(\Omega)$ for all $t\ne0$.
	 	\end{enumerate}	
\end{assumption}
 $\mathcal{R}_{L,6}$ \eqref{it:C1} implies that  $t\mathbb{S}_{\mt{X}}$ is a set of uniqueness of quasi-analytic classes for all $t\ne0$. 

In order to have quantitative statements on $\mathbb{S}_{\boldsymbol{X}}$ yielding sets of uniqueness and explicit conditions so that $\mathcal{R}_{L,6}$ \eqref{it:C2} holds, we restrict our attention to the case where there are sets of points $V_{k}\subseteq\R$ for $k\in[p]$ such that $\prod_{k=1}^pV_{k} \subset \mathbb{S}_{\mt{X}}$. 
The following restriction involves pairs $(C_k,\mathcal{P}_k)$ for $k\in[p]$.  $\mathcal{P}_k$ are classes of probabilities of the form
$$\mathcal{P}_k=\left\{\mathbb{P}_{\boldsymbol{\beta}_k}\in\mathfrak{M}_1^*(\Omega_k): \forall h\in\mathcal{H}_k,\ \mathbb{E}\left[h\left(\boldsymbol{\beta}_k\right)\right]\leq M_k(h)\right\},$$
where $\Omega_k\subseteq\R$ is a closed set, which can be taken not equal to $\R$ if the researcher knows $\mathbb{S}_{\boldsymbol{\beta}_k}\subseteq \Omega_k$, and $\mathcal{H}_k$ can consist of one or a sequence of nonnegative measurable functions. The parameters of these classes are $\Omega_k$, $\mathcal{H}_k$, and $M_k$. 
The proof relies on the related class
$$\mathcal{P}_{c,k}=\left\{\mu\in\mathfrak{M}_c^*(\Omega_k): \forall h\in\mathcal{H}_k,\ \int_{\R}h(z)d|\mu|(z)\leq 2M_k(h)\right\}.$$
The classes $C_k$ are such that $ \mathcal{F}[\mathcal{P}_{c,k}]\subseteq C_k\subseteq C^{\infty}(\R)$. 
\begin{assumption} $\mathbb{S}_{\boldsymbol{X}}$ and $\mathcal{P}$ are such that, for all $k\in[p]$, 
\begin{enumerate}[$\mathcal{R}_{L,7}$\ \textup{(}i\textup{)}]
\item\label{R61} $tV_k$ is a set of uniqueness of $C_k$ for all $t\ne0$; 
\item\label{R62} $\mathbb{P}_{\boldsymbol{\beta}_k}\in\mathcal{P}_k$. 
\end{enumerate}	
\end{assumption}

Let us give examples of pairs $(V_k,\mathcal{P}_k)_{k=1}^p$ and of corresponding sets of uniqueness. Some examples use log-convex sequences $(M_m)_{m\in\N_0}$ which are sequences  of nonnegative numbers such that, for all $m\in\N$, $M_m^2\le M_{m-1}M_{m+1}$ 
 and we also assume $M_0=1$. Define also the iterated logarithms by $\log^{*0}(t) = t$ and for $j\ge1$ by $\log^{*j}(t) = \log(\log^{*(j-1)} (t))$ provided $t$ is large enough for the quantity to be defined and $\exp^{*n}$ is the reciprocal function. 
\begin{enumerate}[\textup{(}{E.}1\textup{)}]
\item\label{e1}   For $k\in [p]$, $\E\left[e^{r\abs{\mt{\beta}_k}}\right] \leq m(r) $ for all $r\in\R_+$, where $ \lim_{r\to \infty}m(r) = \infty $ and $V_{k}$ satisfies 
\begin{equation}\label{eq:e1}
\forall t>0,\ \exists \alpha>1:\ \overline{\lim}_{r\to\infty}\frac{\log(\alpha)}{\log(m(\alpha t r))}\left|V_k\cap \left((-r,r)\setminus\{0\}\right)\right|>1.
\end{equation}
\item\label{e2}  For $k\in [p]$, there exists $(b,c)\in \R_+^2$ such that $\E\left[\abs{\mt{\beta}_k}^m\right] \leq c b^m m! $  for all $m\in \N_0$ and $V_k$ satisfies
\begin{equation}\label{eq:A10}
\overline{\lim}_{r\to\infty}\frac{\log\left(\left|V_k\cap (-r,r)\right|\right)}{r}=\infty.
\end{equation}
\item\label{e3}
 For $k\in [p]$, $\mathbb{S}_{\boldsymbol{\beta}_k}\subseteq \R_+$, there exist $(b,c)\in \R_+^2$ and a log-convex sequence $\{M_m\}$ such that $M_{0}=1$, $\E\left[\abs{\mt{\beta}_k}^m\right] \leq c b^m M_{m} $  for all $m\in \N_0$, 
\begin{equation}\label{it:sign}  \sum_{m\in\N}\frac{1}{M_{m}^{1/(2m)}}=\infty,
\end{equation}	
and $V_{k}$ contains a bounded sequence of distinct points;
\item\label{e4} For $k\in [p]$, there exist $(b,c)\in \R_+^2$ and a log-convex sequence $\{M_m\}$ such that $M_{0}=1$, $\E\left[\abs{\mt{\beta}_k}^m\right] \leq c b^m M_{m} $  for all $m\in \N_0$, and either
	\begin{enumerate}[\textup{(}{E.4}a\textup{)}] 
	\item\label{e4a}  $M_{m}= \nu(m)m!$,  where $\nu$ satisfies, for all $n\in\N$, $\nu(m)=1$ for $0\le m<\exp^{*n}(1)$ and else $\nu(m)=(\log(m)\log^{*2}(m)\cdot\log^{*n}(m))^m$, and $V_k$ satisfies 
	$$\overline{\lim}_{r\to\infty}\frac{\log^{*n+1}\left(\left|V_k\cap (-r,r)\right|\right)}{r}=\infty;$$
		\item\label{e4b} or  $ \{M_{m}\}$ satisfies
	\begin{equation}\label{ideni}
	\sum_{m\in\N}\frac{1}{M_{2m}^{1/(2m)}}=\infty
	\end{equation} 
and $V_{k}$ contains a bounded sequence of distinct points;
\end{enumerate}
\end{enumerate}
In (E.\ref{e1}), a sufficient condition for \eqref{eq:e1}  when $\mathbb{S}_{\boldsymbol{\beta}_k}\subseteq[-\rho,\rho]$, hence $m(r)=e^{\rho r}$, is 
$$\overline{\lim}_{r\to\infty}\frac{1}{r}\left|V_k\cap \left((-r,r)\setminus\{0\}\right)\right|=\infty.$$
In (E.\ref{e2}),  $C_{k}$ consists of the functions $f$ such that there exists $\rho>0$ such that $f$ can be extended uniquely to $\mathcal{H}^{\infty}\left(\{z\in \C:\ |\mathrm{Im}(z)|<\rho\}\right)$ (see Theorem 19.9 in \cite{Rudin2}).

In  (E.\ref{e3}), and (E.\ref{e4}), $C_{k}$ are quasi-analytic classes which are detailed in Appendix \ref{app:complement}. In (E.\ref{e3}), using Proposition \ref{t:PWS2} \eqref{it:Hamb} with $\{M_m\}=\{m!\}$, $C_{k}$ is an analytic class and the set of uniqueness is given in (2) in Section 1 of \cite{Hirschman}. In particular, Proposition \ref{t:PWS2} shows that $ \mathcal{F}[\mathcal{P}_{c,k}]\subseteq C_k$. Based on  Theorem 4b in \cite{Hirschman}, Example 2 in Appendix \ref{app:complement} gives more general forms of $\nu$ in the definition of $\{M_m\}$ in (E.\ref{e4a}) and associated sets $V_k$. The relation between condition \eqref{ideni} (resp. \eqref{it:sign}) and the determinacy of measures in $\mathfrak{M}^*(\R)$ (resp. $\mathfrak{M}^*(\R_+)$) is explicited in Section \ref{sec:CD}. 
\begin{remark}\label{r:CM}
The Student's $t$ with $0<\nu<\infty$ degrees of freedom, the generalized gamma $GG(a,b,p)$ for $0<a<1/2$ of density $ab^px^{ap-1}\exp(-bx^a)/\Gamma(p)$ on $\R^+$, any positive power of the lognormal, the law of $N^{2n+1}$ for all $n\in\N$, $|N|^r$ for $r>4$, $X^m$ for all $m\in\N\setminus\{1,2\}$, $Y^r$ for all $|r|>2$, where $N$ is a Gaussian, $X$ 
a Laplace, gamma or logistic, and $Y$ an inverse Gaussian random variables, respectively do not satisfy these conditions. However, $|N|^r$ for all $0<r\le4$, $X^m$ for $m=1,2$, $Y^r$ for all $-2\le r\le 2$, $GG(a,b,p)$ for all $a\ge1/2$ (thus the $\chi^2$ with any degrees of freedom) satisfy these conditions, hence are determinate in the space of measures with the appropriate support restriction (\emph{i.e.}, $\mathfrak{M}^*(\R)$ or $\mathfrak{M}^*(\R_+)$) (see \cite{Stoyanov,Pakes,Stoy} for more examples).
 \end{remark}

\section{The random coefficients binary choice model}\label{s:BC}
The random coefficients binary choice model takes the form 
\begin{equation}\label{eBCRC0}
Y=\indic\{\alpha+\boldsymbol{\beta}^{\top}\boldsymbol{X}\ge0\}, \ (\alpha,\mt{\beta}) \perp \mt{X}.
\end{equation}
The binary variable $Y$ is 1 when an individual chooses a good, treatment, or an action. 
In the context of treatment effects, $Y=1$ when an individual chooses treatment. \cite{Vytlacil} shows  
the monotonicity of \cite{IA} is equivalent to a selection equation with an additively separable latent index with a single unobservable. Hence monotonicity can be related to rank invariance (see \cite{CH}) similarly to the fact that, in Section \ref{sec:Fourier}, a linear random coefficients models where the random coefficients are functions of a scalar unobservable can be related to the quantile regression. It is well known that monotonicity is sometimes too  strong and it is exemplified in \cite{IA}.  \cite{HV} call \eqref{eBCRC0} the benchmark nonseparable, nonmonotonic model of treatment choice.
Because the scale of the index is not identified, various normalizations can be used. 
A simple one used in consists in assuming one coefficient is 1. This requires that in the original scale the coefficient has a (strict) sign a.s. Based on this and under sufficient variation of the corresponding (special) regressor \cite{Gautier2} show that  identification corresponds to identification of the model of Section \ref{sec:Fourier}. It was mentioned that this idea applies to all sorts of models involving as constitutive element random coefficients indices lying in certain rectangles  (\emph{e.g.}, choice models with multiple alternatives or entry games). \cite{Gautier2} also considers the case where the 
special regressor and the remaining ones have limited variation.  The scaling of \cite{GK} in a preliminary version  \cite{Gautier2} was removed and overlooked as not important. However, assuming that one coefficient is 1 is unnecessary and too restrictive because it still imposes some form of monotonicity. Some elements can be found in \cite{Gautier} in the context of selection equations and missing data in sample surveys. 

Assuming that 
$\mathbb{P}(|(\alpha,\boldsymbol{\beta}^{\top})^{\top}|_2=0)=0$, \eqref{eBCRC0} can be equivalently written as
\begin{equation}\label{eBCRC}
Y=\indic\{\boldsymbol{\Gamma}^{\top}\boldsymbol{S}\ge0\}, \ \mt{\Gamma} \perp \mt{S},
\end{equation}
where $\boldsymbol{\Gamma}=(\alpha,\boldsymbol{\beta}^{\top})^{\top}/|(\alpha,\boldsymbol{\beta}^{\top})^{\top}|_2$ and $\boldsymbol{S}=(1,\boldsymbol{X}^{\top})^{\top}/|(1,\boldsymbol{X}^{\top})^{\top}|_2$. Clearly $|(1,\boldsymbol{X}^{\top})^{\top}|_2\ge1$ and the support of $\boldsymbol{S}$ is a closed subset of the hemisphere $H^+=\{\boldsymbol{s}\in\mathbb{S}^p:\ \boldsymbol{s}_1\ge0\}$. 
We consider identification of the density $f_{\boldsymbol{\Gamma}}^*$ of $\mathbb{P}_{\boldsymbol{\Gamma}}^*$ with respect to $\sigma$, which is the surface measure on $\mathbb{S}^p$.
In this section, we consider the following restriction of the class $\mathcal{P}$
$$\mathcal{P}_{BC}=\left\{\mathbb{P}_{\boldsymbol{\Gamma}}\in\mathfrak{M}_1(\mathbb{S}^p):\ d\mathbb{P}_{\boldsymbol{\Gamma}}=f_{\boldsymbol{\Gamma}}d\sigma,\ 
f_{\boldsymbol{\Gamma}}(\boldsymbol{u})f_{\boldsymbol{\Gamma}}(-\boldsymbol{u})=0 \  \mathrm{for\ a.e.}\ \boldsymbol{u}\in\mathbb{S}^p
\right\}.$$	
It is shown in \cite{GLP} that $\mathbb{P}_{\boldsymbol{\Gamma}}^*$ is identified under the restriction $\mathbb{S}_{\boldsymbol{X}}=\R^{p}$ and $\mathcal{P} =\mathcal{P}_{BC}$.
It is assumed in \cite{GK} that the support of $\boldsymbol{\Gamma}$ lies in an (unknown) hemisphere, namely, that there exists $\boldsymbol{n}$ in $\mathbb{S}^p$ such that $\mathbb{P}(\boldsymbol{n}^\top\boldsymbol{\Gamma}\ge0)=1$.  This assumption first appeared in \cite{IT} and is just a sufficient but unnecessary assumption, similar to assuming that one coefficient is 1. Indeed, if $\boldsymbol{n}\in H^+$, then we have 
$\mathbb{P}(Y=1|\boldsymbol{S}=\boldsymbol{n})=1$, 
else
$\mathbb{P}(Y=1|\boldsymbol{S}=-\boldsymbol{n})=0$. 
This means that there exist limits of values of the regressors such that in the limit everyone chooses $Y=1$ or in the limit everyone chooses $Y=0$. It is stronger than 
$f_{\boldsymbol{\Gamma}}(\boldsymbol{u})f_{\boldsymbol{\Gamma}}(-\boldsymbol{u})=0$ for a.e. $\boldsymbol{u}\ \mathrm{in}\ \mathbb{S}^p$ which does not imply ``unselected samples". 

Consider now the case where the support of $\boldsymbol{X}$ is a proper subset of $\R^p$ or equivalently the support of $\boldsymbol{S}$ is a proper subset of $H^+$. The hemispherical transform is defined, for $\boldsymbol{s}\in\mathbb{S}^p$ and $f\in L^1\left(\mathbb{S}^p\right)$, by 
$\mathcal{T}f(\boldsymbol{s})=\int_{\mathbb{S}^p}\indic\{\boldsymbol{u}^{\top}\boldsymbol{s}\ge0\}f(\boldsymbol{u})d\sigma(\boldsymbol{u})-1/2$. Using the restrictions $\mathcal{P}_{BC}$, we have
\begin{equation}\label{eeg}
\forall \boldsymbol{s}\in\mathbb{S}_{\boldsymbol{S}},\ \mathcal{T}f_{\boldsymbol{\Gamma}}(\boldsymbol{s})+\frac12=\mathbb{E}[Y|\boldsymbol{S}=\boldsymbol{s}].
\end{equation}
$V$-quasi-analytic classes of functions on $\mathbb{S}^p$  are vector spaces of functions on $\mathbb{S}^p$ characterized by $\left(\Delta^m f(\boldsymbol{x})\right)_{m\in\N_0}$ for all $\boldsymbol{x}\in V$. 

In this section, we consider two different identifying restrictions which we denote by $\mathcal{R}_{BC,j}$ for
$j\in [2]$. The main theorem of this section is the following

\begin{theorem}\label{theo:binary}
	$\mathbb{P}_{\boldsymbol{\Gamma}}^*$ in \eqref{eBCRC} is identified under either of $\mathcal{R}_{BC,j}$ for $j\in\{1,2\}$.
\end{theorem}

Similarly to Section \ref{sec:Fourier}, the proofs rely on the diagram \eqref{ediag} with $\mathcal{G}=\mathcal{T}$, $U = \{ \mt{s} = (1,\mt{x}^{\top})/\left|(1,\mt{x}^{\top})\right|_2, \ \mt{x} \in \mathbb{S}_{\mt{X}} \}$, and
$$\Pi : \ \left( \mathbb{P}_{\indic\{\mt{\Gamma}^{\top}\boldsymbol{s}\geq 0\}} \right)_{\mt{s}\in \mathbb{S}_{\mt{S}}} \mapsto \left(\mt{s} \in \mathbb{S}^p \mapsto  \int_{\mathbb{S}^p}\indic\{\boldsymbol{u}^{\top}\boldsymbol{s}\ge0\}f(\boldsymbol{u})d\sigma(\boldsymbol{u})-\frac{1}{2} \right).$$ 
In the rest of this section, we consider $\mathcal{P}\subseteq \mathcal{P}_{BC}$ such that, for given $V\subseteq \mathbb{S}^p$, the vector space spanned by $\mathcal{T}\left[f_{\boldsymbol{\Gamma}}\right]$ for all $\mathbb{P}_{\boldsymbol{\Gamma}}\in\mathcal{P}$ is a $V$-quasi-analytic class of functions on $ \mathbb{S}^p$. We now give details on the restrictions.  

For functions on $\mathbb{S}^p$, the Laplacian $\Delta$ has eigenspaces $H_{m,p}$, eigenvalues $\zeta_{m,p}=-m(m+p-1)$, and $Q_{m,p}f(\boldsymbol{\cdot})=\int_{\mathbb{S}^p}q_{m,p}(\boldsymbol{\cdot},\boldsymbol{y})f(\boldsymbol{y})d\sigma(\boldsymbol{y})$ is the orthogonal projection of $f$ onto $H^{m,p}$ for all $m\in\N_0$. 

\begin{assumption}
	There exists $0<\epsilon<1$ such that 
	\begin{enumerate}[$\mathcal{R}_{BC,1}$\ \textup{(}i\textup{)}]
		\item\label{Rsphere1bb} 
		$\mathbb{S}_{\boldsymbol{X}}$ is a set of uniqueness of $\mathcal{H}^{\infty}\left(\left\{\boldsymbol{z}\in\C^p:\ \left|\mathrm{Im}(\boldsymbol{z})\right|_{2}<\epsilon\right\}\right)$;
		\item\label{Rsphere2bb} $\mathcal{P}$ is 	such that
	$\overline{\lim}_{m\to\infty}\left\|Q_{2m+1,p}f_{\boldsymbol{\Gamma}}\right\|_{L^{1}\left(\mathbb{S}^p\right)}^{1/m}<1/(1+2\epsilon)$.
			\end{enumerate}
\end{assumption}

$\mathcal{R}_{BC,1}$ \eqref{Rsphere2bb}  is a sufficient condition for:  $\E\left[Y | \mt{X} = \cdot\right]$ belongs to  $\mathcal{H}^{\infty}\left(\left\{\boldsymbol{z}\in\C^p:\ \left|\mathrm{Im}(\boldsymbol{z})\right|_{2}<\epsilon\right\}\right)$. Clearly, a set of uniqueness of $\mathcal{H}^{\infty}\left(\left\{\boldsymbol{z}\in\C^p:\ \left|\mathrm{Im}(\boldsymbol{z})\right|_{2}<\epsilon\right\}\right)$ is a set of uniqueness of the superset $\mathcal{H}^{\infty}\left(\left\{\boldsymbol{z}\in\C^p:\ \left|\mathrm{Im}(\boldsymbol{z})\right|_{\infty}<\epsilon\right\}\right)$.
Hence, a sufficient condition for $\mathcal{R}_{BC,1}$ \eqref{Rsphere1bb} is that $U_{p} \subseteq \mathbb{S}_{\boldsymbol{X}}$ where 
for all $j=2,\dots,p$, 
$U_{j}=\bigcup_{\boldsymbol{u}\in{U_{j-1}}}\left\{(\boldsymbol{u}^{\top},v)^{\top},\ v\in {V}_{j}(\boldsymbol{u})\right\}$, where ${V}_{j}(\boldsymbol{u})$ and $U_{1}$ are sets of uniqueness of $\mathcal{H}^{\infty}\left(\left\{\boldsymbol{z}\in\C:\ \left|\mathrm{Im}(\boldsymbol{z})\right|<\epsilon\right\}\right)$. A particular case is a product of sets of uniqueness of $\mathcal{H}^{\infty}\left(\left\{\boldsymbol{z}\in\C:\ \left|\mathrm{Im}(\boldsymbol{z})\right|<\epsilon\right\}\right)$. We give examples in Section \ref{sec:tools}. 

The odd part $f^-$ of $f\in L^q\left(\mathbb{S}^p\right)$ is the limit in $L^q\left(\mathbb{S}^p\right)$ of $f_n^-(\boldsymbol{x})=(f_n(\boldsymbol{x})-f_n(-\boldsymbol{x}))/2$, where $\left(f_n\right)_{n\in \N_0}\in \left(C\left(\mathbb{S}^p\right)\right)^{\N_0}$ converges to $f$ in $L^q\left(\mathbb{S}^p\right)$. We make use of $C^{\infty}_{{\rm odd}}(\mathbb{S}^p)$, the restriction to odd functions of $C^{\infty}(\mathbb{S}^p)$, and
$$C(\{M_m\})=\{f\in C^{\infty}(\mathbb{S}^p):\ \exists c,b\in\R_+:\ \forall m\in\N_0,\ \left\|\Delta^m f^-\right\|_{L^{1}\left(\mathbb{S}^p\right)}\le cb^mM_m\}.$$

\begin{assumption} $\{M_m\}$ satisfies \eqref{ideni} and 
\begin{enumerate}[$\mathcal{R}_{BC,2}$\ \textup{(}i\textup{)}] 
	\item\label{Rsphere1} 
	There exists $U\subseteq\mathbb{S}_{\boldsymbol{S}}$ a set of uniqueness of harmonic homogenous polynomials of odd degree in $\R[Z_1,\dots,Z_{p+1}]$ and, for all $f\in C^{\infty}_{{\rm odd}}\left(\mathbb{S}^p\right)\cap \mathcal{T}[C(\{M_m\})]$, 
if, for all $\boldsymbol{u}\in \mathbb{S}_{\boldsymbol{S}}$,  $f(\boldsymbol{u})=0$, then, for all $\boldsymbol{u}\in U$ and $m\in\N_0$, $\Delta^m f (\boldsymbol{u})=0$;	
	\item\label{Rsphere5} $\mathcal{P}$ is such that 	$ f_{\boldsymbol{\Gamma}} \in C(\{M_m\})$. 
\end{enumerate}
\end{assumption}
$U$ is a set of uniqueness of homogeneous polynomials in $\mathbb{R}[Z_1,\dots,Z_{p+1}]$ if there exists $A\in GL(\mathbb{\R}^{p+1})$ such that
$\{(1,\mt{u}^{\top})^{\top}: \boldsymbol{u}\in \widetilde{U}\}\subseteq AU$ and $\widetilde{U}$ is 
a set of uniqueness of $\mathbb{R}[Z_1,\dots,Z_p]$.  If $\mathbb{S}_{\boldsymbol{X}}$ and thus $\mathbb{S}_{\boldsymbol{S}}$ contains an open set, then clearly $\mathcal{R}_{BC,2}$ \eqref{Rsphere1} is satisfied. 

\section{Panel data models with linear random coefficients model}\label{sec:panel_all}

\subsection{Extension of the Kotlarski lemma}\label{sec:panel}
Consider the equation 
\begin{equation}\label{eq:Kot}
\boldsymbol{Y}_t= \delta + \boldsymbol{\epsilon}_t,\ t\in\{1,2\},
\end{equation}
where the vector of unobserved heterogeneity is $\boldsymbol{\Gamma}:=(\boldsymbol{\epsilon}_1,\boldsymbol{\epsilon}_2,\delta)$. 
\begin{assumption}
	Define the restriction $\mathcal{R}_K$ on $\mathcal{P}$ by $\left\{\mathbb{P}_{\boldsymbol{\Gamma}}=\mathbb{P}_{\boldsymbol{\epsilon}_1}\otimes\mathbb{P}_{\boldsymbol{\epsilon}_2}\otimes\mathbb{P}_{\delta},\ \mathbb{P}_{\boldsymbol{\epsilon}_1}\in \mathcal{P}_1(\Omega_1)\right\}$, where $\Omega_1 \subseteq \R$ is a closed set and $\mathcal{P}_1(\Omega_1)$ is the set of measures which are determinate in $\mathfrak{M}^*(\Omega_1)$ and such that $\E_{\mathbb{P}_{\boldsymbol{\epsilon}_1}}[\boldsymbol{\epsilon}_1]=0$
\end{assumption}

\begin{theorem}\label{Kot}
	$\mathbb{P}_{\boldsymbol{\Gamma}}^*$ in \eqref{eq:Kot} 
	is identified under $\mathcal{R}_K$. 
\end{theorem}

Lemma 1 in \cite{kotlarski} assumes all characteristic functions do not vanish and in Remark 1 it is written that this can be extended to the case where all characteristic functions are analytic. Theorem \ref{Kot} shows that these assumptions are too strong and identification can be achieved when none of the characteristic functions are analytic and $\boldsymbol{\epsilon}_2$ and $\delta$ might not have finite first absolute moments. \cite{evdokimov2012some} present a similar result under alternative assumptions, but assuming $\varphi_{\boldsymbol{\epsilon}_1}$ is analytic.  

\subsection{Linear panel data model where regressors are monomials of a baseline scalar regressor}\label{s:LP}
Consider the equation
\begin{equation}\label{e:panel1}
\boldsymbol{Y}_t=\alpha  + \sum_{j=1}^{T}\boldsymbol{\beta}_{j} \boldsymbol{X}_t^j  +\boldsymbol{\epsilon}_t, \ t=1,\dots,T,
\end{equation}
where $\boldsymbol{X}_t$ is a scalar regressor and denote by $\boldsymbol{\Gamma} = (\alpha,\mt{\beta}^{\top})^{\top}$.  For each 
$t$, we have $$\mathbb{S}_{\boldsymbol{X}_t,\dots,\boldsymbol{X}_t^T}\subseteq\left\{\bold{u}\in\R^p:\ \bold{u}_2=\bold{u}_1^2,\hdots,\bold{u}_T=\bold{u}_1^T\right\},$$
hence the restriction $\mathcal{R}_{L,3}$ does not hold. We show in this section that the availability of $T$ periods allows nonparametric identification. 

\begin{remark}
	Note that \eqref{e:panel1} could be generalized to $\boldsymbol{Y}_t=\alpha  + P(\boldsymbol{X}_t) + \boldsymbol{\epsilon}_t$,  for all $t=1,\dots,T$, where $P(\boldsymbol{X}_t) = \sum_{j=1}^{T}\boldsymbol{\beta}_{j} \boldsymbol{X}_t^{\theta(j)}$ and $\theta \in \N_0^{\N}$ and increasing. However, using our identification strategy would yield to consider the so-called generalized Vandermonde matrices, whose theoretical properties, in particular the ones of their inverse, are not yet well known.
	\end{remark}

\begin{assumption}
	\begin{enumerate}[${\mathcal{R}_{LP,0}}$ \textup{(}i\textup{)}]
			\item\label{R0LPii} $\mathcal{X}_1 =\{\boldsymbol{x}\in\mathbb{S}_{\boldsymbol{X}}: \boldsymbol{x}_1=\dots=\boldsymbol{x}_T\}\neq \emptyset$.
		\item\label{R0LPi} $ \mathbb{P}_{\boldsymbol{\Gamma}} =  \bigotimes_{j=1}^T \mathbb{P}_{\boldsymbol{\epsilon}_j} \otimes \mathbb{P}_{\alpha,\boldsymbol{\beta}}$, $\mathbb{P}_{\boldsymbol{\epsilon}_1}\in \mathcal{P}_1(\Omega_1)$, where $\Omega_1 \subseteq \R$ is a closed set and $\mathcal{P}_1(\Omega_1)$ is the set of measures which are determinate in $\mathfrak{M}^*(\Omega_1)$ and such that $\E_{\mathbb{P}_{\boldsymbol{\epsilon}_1}}[\boldsymbol{\epsilon}_1]=0$;
		\end{enumerate}
\end{assumption}

Using an extension of Theorem \ref{Kot} to $T$ periods, see proof of Theorem \ref{theo:LP2T}, $\mathbb{P}_{\boldsymbol{\epsilon}}^*$ is identified under $\mathcal{R}_{LP,0}$. Note that the restriction ${\mathcal{R}_{LP,0}}$ \eqref{R0LPii} is weaker than assuming that the covariates are centered $\mt{0}_c\in \mathbb{S}_{\boldsymbol{X}}$. The restriction  ${\mathcal{R}_{LP,0}}$ \eqref{R0LPii} is also maintained in \cite{cooprider2020panel}, where they focus on the marginals of $\boldsymbol{\Gamma}$ without imposing our baseline independence assumption but considering only the individuals whose $\mt{X}$ belong to $\mathcal{X}_1$, which are the stayers.

\begin{assumption}
	\begin{enumerate}[${\mathcal{R}_{LP}}$ \textup{(}i\textup{)}]
		\item\label{R1LPiiiT} For all $t\in[T]$, $\mathbb{S}_{\mt{X}_t}$ contains a bounded sequence of distinct points;
		\item\label{R1LPivT}  
		$\mathbb{P}_{\mt{\beta}_T} \in \mathcal{P}_T(\Omega_T)$, where $\Omega_T \subseteq \R$ is a closed set and $\mathcal{P}_T(\Omega_T)$ is the set of measures which are determinate in $\mathfrak{M}^*(\Omega_T)$.  	\end{enumerate}
\end{assumption}

\begin{theorem}\label{theo:LP2T}  $\mathbb{P}^*_{\boldsymbol{\Gamma},\mt{\epsilon}} $ in \eqref{e:panel1} is identified under $\mathcal{R}_{LP,0}$ and $\mathcal{R}_{LP}$. 
\end{theorem}

Like Theorem \ref{tL} under $\mathcal{R}_{L,6}$ or $\mathcal{R}_{L,7}$, Theorem \ref{theo:LP2T} 
makes no assumptions on $\mathbb{P}^*_{\alpha,\mt{\beta}_1,\dots,\mt{\beta}_{T-1}}$. 
\subsection{A single-index panel data model with two periods.}\label{sec:SI}
Consider the equation
\begin{equation}\label{m:SI}
\boldsymbol{Y}_t=f(\boldsymbol{\Gamma}^{\top}\boldsymbol{X}_t) + \boldsymbol{\eta}_t,\ t=1,2,\ f\text{ is increasing.}
\end{equation}

\begin{assumption}
	Either of $\mathcal{R}_{BC,j}$ for $j\in\{1,2\}$ holds and 	
	\begin{enumerate}[${\mathcal{R}_{SI}}$\textup{(}i\textup{)}]
		\item\label{R0SIii}$\{(\boldsymbol{x}_1,\boldsymbol{x}_2)\in\mathbb{S}_{\boldsymbol{X}_1,\boldsymbol{X}_2}: \boldsymbol{x}_1=\boldsymbol{x}_2\}\neq \emptyset$;
		\item \begin{enumerate}
			\item\label{R0SIi}$ \mathbb{P}_{\boldsymbol{\Gamma}, \mt{\eta}} =  \mathbb{P}_{\boldsymbol{\Gamma}} \otimes \mathbb{P}_{\boldsymbol{\eta}_1} \otimes \mathbb{P}_{\boldsymbol{\eta}_2} $,  $\mathbb{P}_{\mt{\eta}_1} \in \mathcal{P}_1(\Omega_1)$, where $\Omega_1 \subseteq \R$ is a closed set and $\mathcal{P}_1(\Omega_1)$ is the set of measures which are determinate in $\mathfrak{M}^*(\Omega_1)$, 
			 $\E_{\mathbb{P}_{\boldsymbol{\eta}_1}}\left[\boldsymbol{\eta}_1\right]=0$; 			\item\label{R0SIiii}$\left|\boldsymbol{\Gamma}\right|_2=1$.
		\end{enumerate}
		
	\end{enumerate}
\end{assumption} 

\begin{theorem}\label{cor:SI}
	$\mathbb{P}_{\boldsymbol{\Gamma},\mt{\eta}}^*$ in \eqref{m:SI} is identified under $\mathcal{R}_{SI}$. \end{theorem}

Theorem \ref{cor:SI}  relies on both theorems \ref{theo:binary} and \ref{Kot}. Restriction ${R_{SI}}$ \eqref{R0SIii} means that there exist ``stayers" in the population, for which the value of the covariate stays the same accross periods, which is a mild assumption. 

%
\renewcommand{\theequation}{A.\arabic{equation}}
\renewcommand{\thelemma}{A.\arabic{lemma}}
\renewcommand{\thecorollary}{A.\arabic{corollary}}
\renewcommand{\thedefinition}{A.\arabic{definition}}
\renewcommand{\theproposition}{A.\arabic{proposition}}
\renewcommand{\theremark}{A.\arabic{remark}}
\renewcommand{\thetheorem}{A.\arabic{theorem}}
\renewcommand{\theassumption}{A.\arabic{assumptio}}
\renewcommand{\thesubsection}{A.\arabic{subsection}}
\setcounter{equation}{0}  
\setcounter{lemma}{0}
\setcounter{corollary}{0}
\setcounter{proposition}{0}
\setcounter{remark}{0}
\setcounter{definition}{0}
\setcounter{lemma}{0}
\setcounter{theorem}{0}
\setcounter{assumption}{0}
\setcounter{section}{0}
\setcounter{subsection}{0}
\setcounter{footnote}{0}
\setcounter{figure}{0}

\section*{\bf APPENDIX}\label{sec:appendix}

\paragraph{\bf Notations.} For a differentiable function $f$ of real variables, $f^{(\boldsymbol{m})}$ denotes $\prod_{j=1}^{p}\frac{\partial^{\boldsymbol{m}_j}}{\partial \mt{x}_j^{\boldsymbol{m}_j}}f$ and $\mathrm{supp}(f)$ its support. Define, for $\mu,\nu\in\mathfrak{M}_c^*\left(\Omega\right)$, where $\Omega$ is a closed set of $\R^p$, the equivalence $\mu\sim\nu$ if $s_{\mu}(\boldsymbol{m})=s_{\nu}(\boldsymbol{m})$ for all $\boldsymbol{m}\in\N_0^{p}$.   $\Aff(\mathbb{R}^p)$ for the group of invertible affine transformations of $\mathbb{R}^p$. Denote by $S'\left(\mathbb{S}^p\right)$ the dual of $C^{\infty}\left(\mathbb{S}^p\right)$ (see \cite{GK}).

\subsection{Complements}\label{app:complement}
\subsubsection{Analytic and quasi-analytic classes of Fourier transforms of complex measures}
Given $b\in\R_+$, and $\{\boldsymbol{M}_{m}\}\in ((0,\infty]^p)^{\N_0}$ where $\{\boldsymbol{M}_{m,k}\}$ for all $k\in[p]$ are log-convex sequences with $\mt{M}_{0,k}=1$,  let us introduce the class of functions
\begin{equation}\label{Cset11}
C^{\{\boldsymbol{M}_{m}\}}(b)=\left\{f\in C^{\infty}(\R^p):\ \exists c:\ \forall \boldsymbol{m}\in\mathbb{N}_0^{p},\ \left\|f^{(\boldsymbol{m})}
\right\|_{L^{\infty}(\R^p)}\le cb^{|\boldsymbol{m}|}\prod_{k=1}^{p}\boldsymbol{M}_{k,\boldsymbol{m}_k}\right\}.
\end{equation}
It consists of complex valued functions and is a vector space. We denote by $C^{\{\mt{M}_{m}\}}=\bigcup_{b\in\R_+}C^{\{\mt{M}_{m}\}}(b)$. Like $\mathcal{A}(S)$, $C^{\{\mt{M}_{m}\}}$ is an algebra with respect to multiplication(see \cite{Rudin2}). 
From Theorem 19.9 in \cite{Rudin2}, $C^{\{m!\}}$ consists of the functions $f$ such that there exists $\rho>0$ such that $f$ can be extended uniquely to $\mathcal{H}^{\infty}\left(\{z\in \C:\ |\mathrm{Im}(z)|<\rho\}\right)$. Note that if $f\in C^{\{\boldsymbol{M}_{m}\}}(b)$ and $\underline{\boldsymbol{x}}\in\R^p$, then the function obtained by fixing all variables but the $k^{\text{th}}$ to those of $\underline{\boldsymbol{x}}$ defines a function of 
$C^{\{\boldsymbol{M}_{k,{\boldsymbol{m}}_k}\}}(b)$. Thus, products of sets of uniqueness of the spaces of functions of a single variable are sets of uniqueness of $C^{\{\boldsymbol{M}_m\}}(b)$.   

By the H\"older inequality, if $\|\cdot\|$ is a norm and $\mu\in\mathfrak{M}^*\left(\R^{p}\right)$, $\left(\int_{\R^{p}}\|\boldsymbol{x}\|^md\mu(\boldsymbol{x})\right)_{m\in\N_0}$ is log-convex and $M_0=\mu\left(\R^{p}\right)$.  Denote by $M(x)$  the trace function of $\{M_m\}$, where, for all $m\in\N$, $M_m^c=\exp\left(\sup_{x\ge0}(mx-M(x))\right)$ and $M(x)=\sup_{m\in\N}(mx-\log(M_m))$ (see \cite{Mand}).  The sequences $\{1_m\}$ and $\{M^c_m\}$ are log-convex. If $S=\R$ and $\underline{\lim}_{m\to\infty}M_m^{1/m}=0$ then $C^{\{M_m\}}(S)=C^{\{0_m\}}(S)$ and if $0<\underline{\lim}_{m\to\infty}M_m^{1/m}<\infty$ then $C^{\{M_m\}}(S)=C^{\{1_m\}}(S)$, else, if the terms in the sequence are positive, $C^{\{M_m\}}(S)=C^{\{M^c_m\}}(S)$.

Let us also introduce certain subsets of $\mathcal{P}_c(\Omega)$ of the form 
$$\mathcal{P}_c^{\{\mt{M}_{m}\}}(\Omega)= 
\left\{  \mu\in\mathfrak{M}_c^*\left(\Omega\right): \ \exists c,b : \forall m\in \N_0, k\in [p],  \  s_{|\cdot|,|\mu|}(me_k) \leq c  b^m\mt{M}_{k,m} \right\}.$$

Proposition \ref{t:PWS2} shows that $\mathcal{F}\left[\mathcal{P}_{c} \right] \subseteq C$ for the classes $\mathcal{P}_{c}$ and $C$ used in the examples (E.\ref{e1}), (E.\ref{e3}) and (E.\ref{e4}) related to the restriction $\mathcal{R}_{L,7}$. 

\begin{proposition}\label{t:PWS2} Let $\Omega \subseteq \R$ be a closed set, $\mu \in \mathfrak{M}_c^*(\Omega)$, $f = \mathcal{F}[\mu]$, and $\{M_m\}\in(0,\infty]^{\N_0}$ be log-convex sequence. 
	\begin{enumerate}
		\item\label{P83b} If there exists an increasing function $m$ such that,   
		for all $r\in\R_+$, $\int_{\R} e^{r|x|}d|\mu|(x)\le 2 m(r)$, then $f$ belongs to 
		$$\mathcal{A}_0(\C)=\left\{f\in\mathcal{A}(\C):\ \forall r\ge0,\ \max_{z:|z|=r}|f(z)|\le 2 m(r)\right\};$$
		\item\label{it:Hamb} Let $\mu\in\mathcal{P}_c^{\{M_{m}\}}(\R)$, where
		$\{M_m\}$ satisfies \eqref{ideni}
		then $f \in C^{\{M_m\}}$, which is a quasi-analytic class.
		\item\label{it:Stiel}  Let $\mu\in\mathcal{P}_c^{\{M_{m}\}}(\R_+)$, where $\{M_m\}$ satisfies \eqref{it:sign}
		then $f \in C^{\{M_m\}}$ and $$C^{\{M_m\}}\bigcap\left\{\mathcal{F}[\mu],\ \mu\in \mathfrak{M}_c^*(\R_+)\right\}$$ is a quasi-analytic class.	
	\end{enumerate}
\end{proposition}

In all the cases of Proposition \ref{t:PWS2}, $f$ is characterized by $\left(f^{m}(0)\right)_{m\in\N_0}$, hence $\mu$ is determinate in $\mathfrak{M}^*_c(\R)$ by injectivity of the Fourier transform.

\begin{proof}[\bf Proof of Proposition \ref{t:PWS2}.]
	Consider \eqref{P83b}. Use that there exists a complex Borel function $g$ with $|g|=1$ such that $d\mu=gd|\mu|$ and rewrite all integrals as integrals with respect to $|\mu|$. Let $r\in \R_+$ and $S_r =\{ z \in \mathbb{C}: \ |z|\leq r\}$.
	For all $x\in\supp(\mu)$, $z\in S_r\mapsto e^{izx}$ is holomorphic and, for all $z\in S_r$, $x\in\R\mapsto e^{izx}$ is $|\mu|$-integrable.
	For all compact $K\subset S_r$, for all $z\in\ K$, we have, for all $x\in\supp(\mu)$, $|e^{izx}|\le e^{r|x|}$ which is $|\mu|-$ integrable. 
	The rest follows by the same argument as those 
	of the proof of theorem p91 of \cite{Tauvel} for complex variables. \\
		Consider \eqref{it:Hamb}. Because $\mu\in\mathfrak{M}_c^*(\R)$, we have $f\in C^{\infty}(\R)$ but also clearly $f\in C^{\{M_m\}}$ because, for all $m\in\N_0$ and $x\in\R$, $\left|f^{(m)}(x)\right|\leq s_{\abs{\cdot},\abs{\mu}}(m)$. 
		The Denjoy-Carleman theorem (see Theorem 19.11 in \cite{Rudin2}) yields that $C^{\{M_m\}}$ is quasi-analytic. Moreover, using that for all $m\in\N_0$, $f^{(m)}(0)=i^{m}s_{\mu}(m)$, $\mu$ is determinate in $\mathfrak{M}_c^*(\R)$.\\ 
		Consider \eqref{it:Stiel}, we adapt the proof of Theorem 4.1 in \cite{CP}. Define the measure $\mu_1$ on $\R$ by $ d\mu_1 (t) =d\mu (t^2)$. Thus, we have, for all $m=2n$ with $n\in\N_0$, $ s_{\abs{\cdot},\abs{\mu_1}}(m) = 2 s_{\abs{\cdot},\abs{\mu}}(n)$ hence $\sum_{m\in\N}1/( s_{\abs{\cdot},\abs{\mu_1}}(2m))^{1/(2m)} = \sum_{m\in\N}1/(2s_{\abs{\cdot},\abs{\mu}}(m))^{1/(2m)}  $ and $\sum_{m\in\N}1/(s_{\abs{\cdot},\abs{\mu}}(m))^{1/(2m)} \geq  \sum_{m\in\N}1/M_{m}^{1/(2m)} = \infty $. Hence, applying \eqref{it:Hamb} to $\mu_1$, $ \mu $ is determinate in $\mathfrak{M}_c^*(\R_+)$. Because we have, for all $ (m,x)\in\N_0\times \R$, 
	\begin{align*}
	\left|\mathcal{F}^{(m)}[\mu](x) \right| & = \df{1}{2} \left|\mathcal{F}^{(2m)}[\mu_1](x) \right| \leq 2^{2m-1} s_{\abs{\cdot},\abs{\mu_1}}(2m)\leq 2^{2m} s_{\abs{\cdot},\abs{\mu}}(m), 
	\end{align*}
	$\mathcal{F}[\mu] \in C^{\{M_m\}}$, hence $\mathcal{F}^{(m)}[\mu](0)=i^{m}s_{\mu}(m)$ yields that $$C^{\{M_m\}}\bigcap\left\{\mathcal{F}[\mu],\ \mu\in \mathfrak{M}_c^*(\R_+)\right\}$$ is quasi-analytic.	
\end{proof}

\subsubsection{Criteria for determinacy}\label{sec:CD}
When $\mu\in\mathfrak{M}^*(\R^{p+1})$, the determinacy of $\mu$ can be assessed by checking the Cram\'er condition:
\begin{equation}\label{e:cramer} \exists \rho>0:\ 
\forall \boldsymbol{t}:\ |\boldsymbol{t}|_2\le \rho,\ \int_{\R^{p+1}}e^{\boldsymbol{t}^{\top}\boldsymbol{x}}d\mu(\boldsymbol{x})<\infty.
\end{equation}
When $\mu\in\mathfrak{M}_1^*(\R^{p+1})$, this means that the moment generating function exists for small $|\boldsymbol{t}|_2$. 
Using that $\int_{\R^{p+1}}\cosh(|{\boldsymbol{t}^{\top}\boldsymbol{x}}|)d\mu(\boldsymbol{x})=\int_{\R^{p+1}}\cosh({\boldsymbol{t}^{\top}\boldsymbol{x}})d\mu(\boldsymbol{x})<\infty$, \eqref{e:cramer} is equivalent to:
\begin{equation}\label{e:cramerb} \exists \rho>0:\ 
\forall \boldsymbol{t}:\ |\boldsymbol{t}|_2\le \rho,\ \int_{\R^{p+1}}e^{|\boldsymbol{t}^{\top}\boldsymbol{x}|}d\mu(\boldsymbol{x})<\infty.
\end{equation}
The Cram\'er condition rules out laws with  heavy tails and holds for example for multivariate normals and skew-normals (see \cite{Stoy} and references therein). For a random variable $X$ such that $\mathbb{S}_{X}\subseteq\R_+$, it is possible to use the weaker Hardy condition $\E[e^{t\sqrt{X}}]<\infty$ for small $|t|$ to assess determinacy (see \cite{Hardy1,Hardy2,StoyanovLin}). We now give a simple extension of Theorem 2 in \cite{Lin} to $\mathfrak{M}^*(\R)$.

\begin{proposition}\label{p:hardy}
	Let $\mu\in\mathfrak{M}^*(\R)$ with $\supp(\mu)=\R_+$, then 
	$\overline{\lim}_{m\to\infty}s_{\mu}(m+1)/(m^2s_{\mu}(m))<\infty$
	implies the following equivalent statements: 
	\begin{enumerate}[\textup{(}{P.A.2.}a\textup{)}] 
		\item\label{H1} $\exists t>0:$ $\int_{\R_+}e^{t\sqrt{x}}d\mu(x)<\infty$;
		\item\label{H2} $\overline{\lim}_{m\to\infty}s_{\mu}(m)^{1/(2m)}/m<\infty$;
		\item\label{H3}  $\overline{\lim}_{m\to\infty}\left(s_{\mu}(m)\right)^{1/m}/(2m)!<\infty$;
	\end{enumerate}
	which imply
	\begin{equation}\label{e:carlemanH}
	\sum_{m\in\N}\frac{1}{s_{\mu}(m)^{\frac{1}{m}}}=\infty,
	\end{equation}
	hence that $\mu$ is determinate. 
\end{proposition}
\eqref{e:carlemanH} is the Carleman condition of the Stieltjes moment problem. 
\begin{corollary}[Theorem 1 in \cite{Lin}]\label{c:carlemanC}
	Let $\mu\in\mathfrak{M}^*(\R)$, then 
	$\overline{\lim}_{m\to\infty}s_{\mu}(2(m+1))/(m^2s_{\mu}(2m))<\infty$
	implies the following equivalent statements:
		\begin{enumerate}[\textup{(}{C.A.1.}a\textup{)}] 
		\item\label{C10} $\exists t\in\R:$ $\int_{\R_+}e^{tx}d\mu(x)<\infty$;
		\item\label{C1} $\exists t>0:$ $\int_{\R_+}e^{t|x|}d\mu(x)<\infty$;
		\item\label{C2} $\overline{\lim}_{m\to\infty}s_{\mu}(2m)^{1/(2m)}/2m<\infty$;
		\item\label{C3} $\overline{\lim}_{m\to\infty}\left(s_{\mu}(2m)/(2m)!\right)^{1/m}<\infty$;
		\item\label{C4} $\overline{\lim}_{m\to\infty}\left(s_{|\cdot|,\mu}(m)/m\right)^{1/m}<\infty$;
	\end{enumerate}
	which imply
	\begin{equation}\label{e:carlemanC}
	\sum_{m\in\N}\frac{1}{s_{\mu}(2m)^{\frac{1}{m}}}=\infty,
	\end{equation}
	hence that $\mu$ is determinate.
\end{corollary}
\eqref{e:carlemanC} is the Carleman condition of the Hamburger moment problem. Corollary \ref{c:carlemanC} (see also  Corollary \ref{c1} below) allows the $\overline{\lim}$ in 
(C.A.1.\ref{C2}) and (C.A.1.\ref{C3})  to be taken over the subsequence of even integers. 
\eqref{e:carlemanC} when the support is known to lie in a half line is weaker than \eqref{e:carlemanH} and nonbinding when $\supp(\mu)$ is compact. 
Neither are necessary for determinacy (see \cite{Heyde} p93, \cite{JS} p113).  

Integral criteria for $\mu\in\mathfrak{M}_1^*(\R)$ with a density $f$ can be easier to check than  \eqref{e:carlemanH} and \eqref{e:carlemanC}. 
Due to Theorem 3  in \cite{Petersen}, they can be used on marginals to assess determinacy of multivariate measures. We have the following results from \cite{Lin2}:
\begin{enumerate}[\textup{(}i\textup{)}]
	\item When $f$ is positive and the Krein condition $\int_{\R}-\log(f(x))/(1+x^2)dx<\infty$ holds, $\mu$ is not determinate; while, if $f$ is also even, differentiable, and  there exists $x_0>0$ such that, for $x\ge x_0>0$ and $x\mapsto f(x)$ decreases to 0, $x\mapsto -xf'(x)/f(x)$ increases to infinity (so-called Lin conditions), and $\int_{\R}-\log(f(x))/(1+x^2)dx=\infty$ then \eqref{e:carlemanH} holds.	
	\item When $\supp(\mu)\subseteq\R^+$, $f(x)$ is positive on $[x_0,\infty)$ for $x_0>0$, and the Krein condition $\int_{x_0}^{\infty}-\log(f(x^2))/(1+x^2)dx<\infty$ holds, $\mu$ is not determinate (see \cite{Pakes} and \cite{Pedersen} for weaker conditions on the interval where $f$ is positive); while, if $f$ is positive and differentiable on $\R_+$, and, for $x\ge x_0$, $x\mapsto f(x)$ decreases to 0 and $x\mapsto-xf'(x)/f(x)$ increases to infinity (so-called Lin conditions), and $\int_{\R_+}-\log(f(x^2))/(1+x^2)dx=\infty$, then \eqref{e:carlemanC} holds.
	\end{enumerate}
The Krein condition is  not necessary for $\mu\in\mathfrak{M}_1^*(\R)$ to be indeterminate (see \cite{JS} p114) in the Hamburger and Stieltjes cases.

Theorem 3 in \cite{Petersen} and Proposition \ref{t:PWS2} imply the following results in higher dimensions.
\begin{corollary}\label{c1}
	Let $q\in\{0,1,\dots,p+1\}$, $\mu_1,\mu_2\in\mathfrak{M}_c^*(\R^{p+1})$ such that for known $\mathfrak{B}\in GL(\R^{p+1})$ $\supp(\mathfrak{B}_*\mu_1)=\supp(\mathfrak{B}_*\mu_2)=\R^q\times\R_+^{p+1-q}$ and $\mu_1\sim\mu_2$.\\
	We have $\mu_1=\mu_2$ if for $\mathfrak{A}\in\Aff(\R^{p+1})$ leaving invariant $\mathrm{span}\{\boldsymbol{e}(q+1),\dots,\boldsymbol{e}(p+1)\}$
	\begin{align*}
	&\min_{j=1,\dots,q}\sum_{m\in\N_0}\frac{1}{\left(s_{|(\mathfrak{A}\circ\mathfrak{B})_*\mu_1|}(2m\boldsymbol{e}(j))+s_{|(\mathfrak{A}\circ\mathfrak{B})_*\mu_2|}(2m\boldsymbol{e}(j))\right)^{\frac{1}{2m}}}=\infty;\\
	&\min_{j=q+1,\dots,p+1}\sum_{m\in\N_0}\frac{1}{\left(s_{|(\mathfrak{A}\circ\mathfrak{B})_*\mu_1|}(m\boldsymbol{e}(j))+s_{|(\mathfrak{A}\circ\mathfrak{B})_*\mu_2|}(m\boldsymbol{e}(j))\right)^{\frac{1}{2m}}}=\infty;
	\end{align*} 
	or, when $\mu_1,\mu_2\in\mathfrak{M}^*(\R^{p+1})$, 
	\begin{align}
	&\min_{j=1,\dots,q}\sum_{m\in\N_0}\frac{1}{s_{(\mathfrak{A}\circ\mathfrak{B})_*\mu_1}(2m\boldsymbol{e}(j))^{\frac{1}{2m}}}=\infty;\label{e:Carleman1}\\
	&\min_{j=q+1,\dots,p+1}\sum_{m\in\N_0}\frac{1}{s_{(\mathfrak{A}\circ\mathfrak{B})_*\mu_1}(m\boldsymbol{e}(j))^{\frac{1}{2m}}}=\infty\label{e:Carleman2};
	\end{align} 
	where, in both cases, if $q=0$, we remove the first minimum, while, if $q=p$, we remove the second minimum.
\end{corollary}
\noindent {\bf Proof of Corollary \ref{c1}.} When $q=p+1$, the statement is obtained by applying the Denjoy-Carleman theorem to $\mathcal{F}[(\mathfrak{A}\circ\mathfrak{B})_*(\mu_1-\mu_2)]$ (see Theorem 3.2 in \cite{CP} for complex measures and else Theorem 2.3 in \cite{DeJeu1}). When $q<p$, we apply this result to  $(\mathfrak{C}\circ\mathfrak{B})_*\mu_j$ for $j=1,2$, where $\mathfrak{C}((\boldsymbol{x}_1,\dots,\boldsymbol{x}_{p+1})^{\top})=(\boldsymbol{x}_1,\dots,\boldsymbol{x}_q,\sqrt{\boldsymbol{x}_{q+1}},\dots,\sqrt{\boldsymbol{x}_{p+1}})$ and proceed like in the proof of Theorem 4.1 in \cite{CP}.\hfill $\square$\vspace{0.3cm}

The Denjoy-Carleman theorem applied to $\mathcal{F}\left[(\mathfrak{A}\circ\mathfrak{B})_*\mu\right]$ where $\mu\in \mathfrak{M}_c\left(\R^{p+1}\right)$ only requires assumptions on marginals, hence the conditions in Corollary \ref{c1} only involve these marginals, this is in line with Theorem 3 in \cite{Petersen}. \cite{Nussbaum} (see also Theorem 5 in \cite{Berg}) states this result for elements of $\mathfrak{M}^*(\R^{p+1})$ when $q=p+1$ without allowing for affine transformations. Due to Proposition 3.10 in \cite{DeJeu2}, it is possible that the sum of two functions belonging to classes defined by two different $p+1$-tuples of sequences which both satisfy the sufficient condition of Denjoy-Carleman theorem  do not belong to such a class. Rather than \eqref{e:Carleman1}-\eqref{e:Carleman2}, also called Carleman conditions, \cite{ST} assumes the stronger condition
\begin{equation}\label{e:carlmauvais}
\sum_{m\in\N_0}\frac{1}{\left(\int_{\R^{p+1}}\sum_{j=1}^{p+1}\boldsymbol{x}_j^{2m}d\mu(\boldsymbol{x})\right)^{\frac{1}{2m}}}=\infty.
\end{equation} 

\begin{corollary}\label{c2} Let $q\in\{0,1,\dots,p+1\}$,  $\mu_1,\mu_2\in\mathfrak{M}_c^*(\R^{p+1})$ such that for known $\mathfrak{B}\in GL(\R^{p+1})$ $\supp(\mathfrak{B}_*\mu_1)=\supp(\mathfrak{B}_*\mu_2)=\R^q\times\R_+^{p+1-q}$ and $\mu_1\sim\mu_2$.\\ 
	$\mu_1=\mu_2$ if there exist functions $W_j$ for $j=1,\dots,p+1$ from $\R$ to $[0,\infty]$ such that:
	\begin{align}
	\min\left\{\min_{j=1,\dots,q}\sum_{m\in\N}\frac{1}{\left\|x\mapsto\frac{x^{2m}}{W_j\left(\sqrt{x}\right)}\right\|_{L^{\infty}(\R)}^{\frac{1}{2m}}},\min_{j=q+1,\dots,p+1}\sum_{m\in\N}\frac{1}{\left\|x\mapsto\frac{x^{2m}}{W_j\left(\sqrt{x}\right)}\right\|_{L^{\infty}(\R_+)}^{\frac{1}{2m}}}\right\}=\infty;\label{cond}
	\end{align}
	where, if $q=0$, we remove the first minimum, while, if $q=p$, we remove the second minimum;  
	and $\mathfrak{A}\in\Aff(\R^{p+1})$ leaving invariant $\mathrm{span}\{\boldsymbol{e}(q+1),\dots,\boldsymbol{e}(p+1)\}$ such that:
	\begin{align*}
	&\sum_{j=1,\dots,p+1}\int_{\R}W_j(\boldsymbol{x}_j)d\left(\Pi[\boldsymbol{e}(j)]_*\left(|(\mathfrak{A}\circ\mathfrak{B})_*\mu_1|+|(\mathfrak{A}\circ\mathfrak{B})_*\mu_2|\right)\right)(\boldsymbol{x})<\infty;\\
	&\mathrm{or}\ \sum_{j=1,\dots,p+1}\int_{\R}W_j(\boldsymbol{x}_j)d\left(\left(\Pi[\boldsymbol{e}(j)]\circ\mathfrak{A}\circ\mathfrak{B}\right)_*\mu_1\right)(\boldsymbol{x})<\infty\ \mathrm{if}\ \mu_1,\mu_2\in\mathfrak{M}^*(\R^{p+1}).
	\end{align*}
	Moreover, $\mu_1=\mu_2$ if there exist functions $W_j$ for $j=1,\dots,p+1$ from $\R$ to $[c_j,\infty]$, where $c_j>0$, satisfying \eqref{cond}, $\mathfrak{A}\in\Aff(\R^{p+1})$ leaving invariant $\mathrm{span}\{\boldsymbol{e}(q+1),\dots,\boldsymbol{e}(p+1)\}$, and a function $W$ from $\R^{p+1}$ to $[0,\infty]$ such that, for $|\boldsymbol{x}|_2$ large enough, 
	$W(\mathfrak{A}\circ\mathfrak{B}(\boldsymbol{x}))\ge\prod_{j=1}^{p+1}W_j(\boldsymbol{x}_j)$,
	with
	$\int_{\R^{p+1}}W(\mathfrak{A}\circ\mathfrak{B}(\boldsymbol{x}))d(|\mu_1|+|\mu_2|)(\boldsymbol{x})<\infty$ or $\int_{\R^{p+1}}W(\mathfrak{A}\circ\mathfrak{B}(\boldsymbol{x}))d\mu_1(\boldsymbol{x})<\infty$ if $\mu_1,\mu_2\in\mathfrak{M}^*(\R^{p+1})$.
\end{corollary}
\noindent {\bf Proof of Corollary \ref{c2}.} Let us prove the moreover part of the statement.  
The first part is obtained in a similar manner. Assuming that $W_j$ are bounded from below by $c_j>0$ for $j=1,\dots,p+1$, we use that for $j\in\{1,\dots,p+1\}$, $m\in\N$, and $\mu,\nu\in\mathfrak{M}_c^*(\R^{p+1})$, we have:
\begin{align*}
s_{|\mu|}(m\boldsymbol{e}(j))+s_{|\nu|}(m\boldsymbol{e}(j))&\le\left\|\boldsymbol{x}^{m}/W_j(\boldsymbol{x})\right\|_{L^{\infty}(\R)}\int_{\R^{p+1}}W_j(\boldsymbol{x}_j)d\left(|\mu|+|\nu|\right)(\boldsymbol{x})\\
&\le \frac{1}{\prod_{\substack{k=1,\dots,p+1\\k\ne j}}c_k}\left\|\boldsymbol{x}^{m}W_j(\boldsymbol{x})^{-1}\right\|_{L^{\infty}(\R)}\int_{\R^{p+1}}W(\boldsymbol{x})d\left(|\mu|+|\nu|\right)(\boldsymbol{x}).
\end{align*}
We use the same argument for $\mu_1\in\mathfrak{M}^*(\R^{p+1})$. We conclude using Corollary \ref{c1}.\hfill $\square$\vspace{0.3cm}

The first sufficient condition in Corollary \ref{c2} is the less demanding. Due to Theorem 3.14 in \cite{DeJeu2}, if $W$ is a function from $\R$ to $[c,\infty]$ for $c>0$, the next statements are equivalent:
\begin{enumerate}[\textup{(}{QAW.}1\textup{)}]
	\item\label{QA0} $\sum_{m\in\N}1/\left\|x\mapsto x^{2m}/W(x)\right\|_{L^{\infty}(\R)}^{1/(2m)}=\infty$;
	\item\label{QA1} There exists a function $\tilde{W}$ with values in $[0,\infty)$ and $R>0$ such that, for $|x|>R$, $W(x)\ge  \tilde{W}(x)$ and $\tilde{W}(x)=\tilde{W}(-x)<\infty$, $x\mapsto \log\left(\tilde{W}(e^x)\right)$ is convex on $(\log(R),\infty)$ and $\int_R^{\infty}\log \tilde{W}(s)/(1+s^2)ds=\infty$;
	\item\label{QA2}  There exist 
	$\epsilon>0$, and a nondecreasing nonnegative function $\rho\in C^{\infty}([0,\infty))$ equal to 0 on $[0,\epsilon]$, such that, for all $x\in\R$, $W(x)\ge \exp\left(\int_0^{|x|}(\rho(s)/s) ds\right)$ and $\int_0^{\infty}\rho(s)/(1+s^2)ds=\infty$.
\end{enumerate} 
For example, if there exist $R\in\R$, $(a_j)_{j\in\N_0}\in(0,\infty)^{\N_0}$, and $(p_j)_{j\in\N_0}\in(-\infty,1]^{\N_0}$ equal to 0 for $j$ large enough, such that
$$W(x)\ge \exp\left(\frac{x^2}{\prod_{j=0}^{\infty}\log_j^{p_j}(a_j|x|)}\right)\indic\left\{|x|\ge R\right\}\quad \left(\mathrm{\emph{e.g.},}\ W(x)=\exp{\left(\frac{|x|}{a_0\log (a_1|x|)}\right)}\indic\left\{|x|\ge R\right\}\right)$$
then $W$ satisfies (QAW.\ref{QA2}).
Replacing $x$ by $\sqrt{x}$ and using that the above lower bound is even, we obtain that $\sum_{m\in\N}1/\left\|x\mapsto x^{2m}/W(\sqrt{x})\right\|_{L^{\infty}(\R_+)}^{1/(2m)}=\infty$ if
there exist $(C,R)\in\R^2$, $(a_j)_{j\in\N_0}\in(0,\infty)^{\N_0}$, and $(p_j)_{j\in\N_0}\in(-\infty,1]^{\N_0}$ equal to 0 for $j$ large enough, such that
$$W(x)\ge C\exp\left(\frac{x}{\prod_{j=0}^{\infty}\log_j^{p_j}(a_j\sqrt{x})}\right)\indic\left\{x\ge R\right\}\quad \left(\mathrm{\emph{e.g.},}\ W(x)=\exp{\left(\frac{\sqrt{x}}{a_0\log (a_1x)}\right)}\indic\left\{x\ge R\right\}\right).$$
Because $\sum_{m\in\N}1/\left\|x\mapsto x^{2m}/W(x)\right\|_{L^{\infty}(\R)}^{1/(2m)}=\infty$ implies $\sum_{m\in\N}1/\left\|\boldsymbol{x}\mapsto(\boldsymbol{x}^{\top}\boldsymbol{f}(j))^{2m}/W(|\boldsymbol{x}|_2)\right\|_{L^{\infty}(\R^q)}^{1/(2m)}=\infty$ for all bases $(\boldsymbol{f}(j))_{j=1}^q$ of $\R^q$ and similarly for $\R_+^{p+1-q}$, assuming that $\mathfrak{B}$ is identity, we can replace \eqref{cond} by the stronger condition (not involving $\mathfrak{A}$)
$$W(\boldsymbol{x})\ge \widetilde{W}\left(|(\boldsymbol{x}_1,\dots,\boldsymbol{x}_q)^{\top}|_2\right)\widetilde{W}\left(\sqrt{|(\boldsymbol{x}_{q+1},\dots,\boldsymbol{x}_{p+1})^{\top}|_2}\right),$$ where $\widetilde{W}$ is a function from $\R$ to $[c,\infty]$ for $c>0$ which satisfies the equivalent (QAW.\ref{QA0})-(QAW.\ref{QA2}).
The function $1/W$ is called a quasi-analytic weight in \cite{DeJeu1,DeJeu2}.

The approach using quasi-analytic weights is closely related to the Krein and Lin conditions. Assume that $\mathfrak{A}$ and $\mathfrak{B}$ are the identity, $q=p+1$ (the remaining cases can be treated similarly), and $\mu\in\mathfrak{M}^*(\R^{p+1})$ has marginals $\Pi[\boldsymbol{e}(j)]_*\mu$ with densities $f_j(\boldsymbol{x}_j)$ for $j=1,\dots,p+1$ which are even and positive and such that $-\boldsymbol{x}_jf_j'(\boldsymbol{x}_j)/f_j(\boldsymbol{x}_j)-\alpha_j\log(\log(\boldsymbol{x}_j))$ with $\alpha_j>1$ is increasing for $|\boldsymbol{x}_j|\ge x_{0j}\ge0$. Then, for $j=1,\dots,p+1$, $W_j(\boldsymbol{x}_j)=1/\max(f_j(\boldsymbol{x}_j)\log(|\boldsymbol{x}_j|)^{\alpha_j},1)$ satisfies (QAW.\ref{QA1}) and $\int_{\R}W_j(\boldsymbol{x}_j)f_j(\boldsymbol{x}_j)d\boldsymbol{x}_j<\infty$, hence the conclusions of Corollary \ref{c2} hold.
The advantage of the approach using quasi-analytic weights is that, unlike the Krein and Lin conditions, it is not required that a density exists.

\subsubsection{Complements on the sets of uniqueness}\label{sec:tools}
Examples of sets of uniqueness for analytic classes are given in the examples (E.\ref{e1}) and (E.\ref{e2}). We only give here more details about (E.\ref{e1}). 
Let $m \in \R_+^{[0,\infty)}$, $\lim_{r\to\infty}m(r)= \infty$, and
$$\mathcal{A}_0(\C)=\left\{f\in\mathcal{A}(\C):\ \exists C>0:\ \forall r\ge0,\ \max_{z:|z|=r}|f(z)|\le C m(r)\right\}.$$
$U$ is a set of uniqueness of such a class $\mathcal{A}_0(\C)$
if, for $\epsilon>0$, we have
\begin{align}\label{eq:Jensen0b}
\exists \alpha>1:\ \overline{\lim}_{r\to\infty}\frac{\log(\alpha)}{\log(m(\alpha r))}\left|U\cap \left((-r,r)\setminus\{0\}\right)\right|>1.
\end{align}
Hence, 
$t U$ is a set of uniqueness for all $t\ne0$ if 
$$\forall t>0,\ \exists \alpha>1:\ \overline{\lim}_{r\to\infty}\frac{\log(\alpha)}{\log(m(\alpha t r))}\left|U\cap \left((-r,r)\setminus\{0\}\right)\right|>1.$$
A sufficient condition when $m(r)=e^{\rho r}$, which occurs when $f$ is the Fourier transform of the difference of two complex measures supported in $[-\rho,\rho]$, is 
$$\overline{\lim}_{r\to\infty}\frac{1}{r}\left|U\cap \left((-r,r)\setminus\{0\}\right)\right|=\infty.$$
Let us explain the above characterization. If $f\in\mathcal{A}_0(\mathbb{C})$, 
there exists $k\in\N$ such that $f(z)=z^kg(z)$ and  $g(0)\ne 0$. Jensen's formula applied to $g$ (see 15.18 in \cite{Rudin2}), with $ \max_{z:|z|=r}|g(z)|= \max_{z:|z|=r}|f(z)|/z^k \leq Cm(r)/r^k $, ensures 
$$  \abs{g(0)}\prod_{k=1}^{n_g(\alpha r)}\dfrac{\alpha r}{|\omega_k|} \leq \df{C m(\alpha r)}{\alpha^k r^k},$$
where $\omega_1,\dots,\omega_{n_g(\alpha r)}$ are the zeros of $g$ in $B(0,\alpha r)$ ranked according to their multiplicity and  $n_g(\alpha r)$ is the number of zeros of $g$ with multiplicity in $B(0,\alpha r)$. Thus, we have 
$$ \df{C m(\alpha r)}{\alpha^k r^k} \geq \abs{g(0)}\prod_{k=1}^{n_g( r)}\dfrac{\alpha r}{|\omega_k|} \geq \abs{g(0)} \alpha^{n_g(r)},$$
hence 
$n_g( r) \leq 
\log(C m(\alpha r)(\alpha r)^{-k}/\abs{g(0)})/\log(\alpha)$.
We have
$$  \df{n_g( r) \log(\alpha)}{\log(m(\alpha r))} \leq  1 + \df{\log(C (\alpha r)^{-k}/\abs{g(0)})}{\log(m(\alpha r))}   $$
then $ \overline{\lim}_{r\to\infty}n_g( r) \log(\alpha)/\log(m(\alpha r)) \leq  1$.  This yields that $U$ satisfying \eqref{eq:Jensen0b} is a set of uniqueness of $\mathcal{A}_0(\C)$.

The following sets $U$ are sets of uniqueness of the following quasi-analytic and $V$-quasi-analytic classes and complement (E.\ref{e4}). \\ 
{\bf Example 1.}  $V $ is a subset of the interior of $U$, which is nonempty. Indeed, a function which is zero on $U$ has all its partial derivatives or Laplacians (for functions on $S\subseteq\mathbb{S}^p$) equal to zero at every point in $V$.\\
{\bf Example 2.}  (Theorem 4b in \cite{Hirschman}) 
Let $M_m= \nu(m)^mm!$, where $x\in[0,\infty)\mapsto\nu(x)$ is increasing and continuously differentiable such that $\nu(0)=1$ and $$\lim_{x\to\infty}\frac{x\nu'(x)}{\nu(x)}=0.$$ 
$U$ is a set of uniqueness of  $C^{\{M_m\}}(b)$ (which contains non-analytic functions) if
\begin{equation}\label{eq:H}
\overline{\lim}_{r\to \infty}\frac1r\int_1^{\left|U\cap (-r,r)\right|}M(\log(u))\frac{du}{u^2}> \frac{\pi b}{2},
\end{equation}
where $M(x)$ is the \emph{trace function} of $\{M_m\}$ defined by $M(x)=\sup_{m\in\N}(mx-\log(M_m))$ (see \cite{Mand}).

\subsubsection{$\mathcal{A}(S)$ is a small subset of $C^{\infty}(S)$}\label{sec:meager}
Proposition \ref{prop:small} can be found in \cite{SV}, it is of the same spirit as the statement that the complement of nowhere analytic functions in $C^{\infty}(S)$ is meager with empty interior, which we detail below. $C^{\{M_m\}}(S)$ is defined like $C^{\{M_m\}}$ but we replacing $\R^p$ by an open set $S$. 
\begin{proposition}\label{prop:small}
	If condition \eqref{ideni} fails, then $C^{\{M_m\}}(S)$ contains a vector space of dimension $2^{\aleph_0}$ such that for all nonzero element $f$ there does not exists $x\in S$, a neighborhood $\mathcal{N}_x$ of $x$, and a sequence $\{M_{x,m}\}$ such that $f$ belongs to a quasi-analytic class $C^{\{M_{x,m}\}}(\mathcal{N}_x)$.	
\end{proposition}

Assume, for simplicity, that $S\subset\R^p$ is compact and as usual $C^{\infty}(S)$ is equipped with the distance 
$d(f,g)=\sum_{\boldsymbol{m}\in\mathbb{N}_0^{p}}\min\left(1/\boldsymbol{2}_c^{\boldsymbol{m}},\left\|f^{(\boldsymbol{m})}-g^{(\boldsymbol{m})}\right\|_{L^{\infty}(S)}\right)$.
By the Cauchy-Hadamard theorem (see \cite{shabat}), $f\in C^{\infty}(S)$ is analytic at $\underline{\boldsymbol{x}}$ if and only if there exist $b,c\in\mathbb{N}$ such that
$f$ belongs to
$$T(\underline{\boldsymbol{x}},b,c)=\left\{f\in C^{\infty}(S):\ \forall \boldsymbol{m}\in\mathbb{N}_0^{p},\ \left|f^{(\boldsymbol{m})}
(\underline{\boldsymbol{x}})\right|\le cb^{|\boldsymbol{m}|}\boldsymbol{m}!\right\},$$
thus 
$\mathcal{A}(S)=\bigcup_{\underline{\boldsymbol{x}}\in\mathbb{Q}^{p}\cap S,b,c\in\mathbb{N}}T(\underline{\boldsymbol{x}},b,c)$.
This is a countable union of closed sets, hence closed. Moreover, $T(\underline{\boldsymbol{x}},b,c)$ has empty interior. Indeed (see, \cite{SZ}), 
given $f$ in $T(\underline{\boldsymbol{x}},b,c)$, for all $\epsilon>0$, $m\in\N$ such that $\sum_{j=m}^{\infty}2^{-j}<\epsilon$, $c<(\epsilon b^{m}/(2m!))^{1/(2m)}$, the functions $ f_{\epsilon}: \ \boldsymbol{x}\in S \mapsto f(\boldsymbol{x})+\epsilon\cos\left(c(\boldsymbol{x}_{1}-\underline{\boldsymbol{x}}_{1})\right)/(2c^m)$
are such that $d(f_{\epsilon},f)<\epsilon$ and $\left|f_{\epsilon}^{(2m)}(\underline{\boldsymbol{x}})-f^{(2m)}(\underline{\boldsymbol{x}})\right|>b^{2m}(2m)!$, hence  $f_{\epsilon}\notin T(\underline{\boldsymbol{x}},b,c)$. Due to Baire's theorem, 
the meager set $\mathcal{A}(S)$ of the complete metric space $C^{\infty}(S)$ has an empty interior. With the arguments in \cite{Cater}, the complement of nowhere analytic functions in $C^{\infty}(S)$, hence containing $\mathcal{A}(S)$, can be shown to be meager with empty interior. 

\subsection{Proofs}\label{sec:proofs}

\begin{proof}[\bf Proof of Theorem \ref{tL} under $\mathcal{R}_{L,1}$.] 
	Assume that $\mathbb{P}_{\boldsymbol{\Gamma}}$ and $\mathbb{P}_{\boldsymbol{\Gamma}}^*$ both give rise to the same collection	$\left(\mathbb{P}_{
	\left(1, \mt{x}^{\top}\right)\boldsymbol{\Gamma}}
	\right)_{\mt{x}\in \mathbb{S}_{\boldsymbol{X}}}$. 
	By \eqref{RL1ii}, we have for all $\boldsymbol{x}\in\mathbb{S}_{\boldsymbol{X}}$ and $t$ in $\R$,
	\begin{equation}\label{eq:prod}
	\varphi_{\alpha}(t)\prod_{j=1}^p\varphi_{\boldsymbol{\beta}_j}(t\boldsymbol{x}_j)=\varphi_{\alpha}^*(t)\prod_{j=1}^p\varphi_{\boldsymbol{\beta}_j}^*(t\boldsymbol{x}_j).
	\end{equation}
	Because the value at 0 of a characteristic function is 1, taking $\boldsymbol{x}=\boldsymbol{0}_c$, yields $\varphi_{\alpha}=\varphi_{\alpha}^*$.\\ 
	Because $\varphi_{\alpha}^*$ is continuous at 0, for all $\boldsymbol{x}\in\mathbb{S}_{\boldsymbol{X}}$ and $t\in(-t_0,t_0)$ for $t_0$ small enough, 
	\begin{equation}\label{eq:cf}\prod_{j=1}^p\varphi_{\boldsymbol{\beta}_j}(t\boldsymbol{x}_j)=\prod_{j=1}^p\varphi_{\boldsymbol{\beta}_j}^*(t\boldsymbol{x}_j).
	\end{equation} 
	Taking $\boldsymbol{x}=\boldsymbol{x}(1)$ yields, for all $t\in(-t_0,t_0)$, 
	$\varphi_{\boldsymbol{\beta}_1}(\boldsymbol{x}(1)_1t)=\varphi_{\boldsymbol{\beta}_1}^*(\boldsymbol{x}(1)_1t)$, 
	hence $\varphi_{\boldsymbol{\beta}_1}(t)=\varphi_{\boldsymbol{\beta}_1}^*(t)$ for all $t\in(-t_0/\boldsymbol{x}(1)_1,t_0/\boldsymbol{x}(1)_1)$. Hence $\mathbb{P}_{\boldsymbol{\beta}_1}$ and $\mathbb{P}_{\boldsymbol{\beta}_1}^*$ have same moments, so $\mathbb{P}_{\boldsymbol{\beta}_1}=\mathbb{P}_{\boldsymbol{\beta}_1}^*$. We conclude by iterating this procedure.\end{proof}
	
	\begin{proof}[\bf Proof of Theorem \ref{tL} under $\mathcal{R}_{L,2}$.] 
	By \eqref{eq:prod}, \eqref{eq:cf} holds for all $\boldsymbol{x}\in\mathbb{S}_{\boldsymbol{X}}$ and $t$ in the complement of the zeros of $\varphi_{\alpha}^*$. Because the complement of the set of zeros of $\varphi_{\alpha}^*$ is dense and $t\mapsto \prod_{j=1}^p\varphi_{\boldsymbol{\beta}_j}(t\boldsymbol{x}_j)$ and $t\mapsto \prod_{j=1}^p\varphi_{\boldsymbol{\beta}_j}^*(t\boldsymbol{x}_j)$ are continuous, \eqref{eq:cf} holds for all $\boldsymbol{x}\in\mathbb{S}_{\boldsymbol{X}}$ and $t\in\R$. Taking $\boldsymbol{x}=\boldsymbol{x}(1)$ yields, for all $t\in\R$, 
	$\varphi_{\boldsymbol{\beta}_1}(\boldsymbol{x}(1)_1t)=\varphi_{\boldsymbol{\beta}_1}^*(\boldsymbol{x}(1)_1t)$, 
	so $\varphi_{\boldsymbol{\beta}_1}=\varphi_{\boldsymbol{\beta}_1}^*$. So, for all $\boldsymbol{x}\in\mathbb{S}_{\boldsymbol{X}}$ and $t\in\R$, 
	$\prod_{j=2}^p\varphi_{\boldsymbol{\beta}_j}(t\boldsymbol{x}_j)=\prod_{j=2}^p\varphi_{\boldsymbol{\beta}_j}^*(t\boldsymbol{x}_j)$, and we conclude by iteration. 
\end{proof}

\begin{proof}[\bf Proof of Theorem \ref{tL} under $\mathcal{R}_{L,3}$.] 
Take $S = \R^{p+1}$ and $\mathfrak{F}(S)$ the vector space spanned by $\mathcal{F}[\mathcal{P}]$.  
An element of $\mathfrak{F}(S)$ is of the form $\mathcal{F}[\mu-\nu]$, where $\mu,\nu\in\mathcal{P}$. %
For all  $k\in \N_0$ and $\mt{x}\in \mathbb{S}_{\mt{X}}$,
\begin{align}\label{eq:moment}
(-i)^k\left. \partial_t^{(k)}\mathcal{F}\left[ \mu-\nu\right](t, t\mt{x})\right|_{t=0} &= \int_{\R^{p+1}} \left(\boldsymbol{\gamma}_1 + \sum_{j=1}^{p}\boldsymbol{\gamma}_{j+1}\mt{x}_j\right)^k d(\mu-\nu)(\boldsymbol{\gamma})
 = P_k(\mt{x}), \notag 
\end{align}
where $P_k$ are the polynomials with real coefficients 
\begin{equation*}
P_k(Z_1,\dots,Z_{p}) = \sum_{\mt{j}\in \N_0^{p+1}: \ |\mt{j}|_1 = k}  \left( \begin{array}{c}
k\\ 
\mt{j}_1,\dots,\mt{j}_{p+1}
\end{array} \right) \left(\int_{\R^{p+1}} \boldsymbol{\gamma}^{\mt{j}}  d(\mu-\nu)(\boldsymbol{\gamma})\right)\  Z_1^{\mt{j}_2}\dots Z_{p}^{\mt{j}_{p+1}}.
\end{equation*}
So, by (i),  $\mathcal{F}\left[\mu-\nu\right]=0$ on $U = \{t(1,\mt{x}): t\in \R, \mt{x}\in \mathbb{S}_{\mt{X}}\}$ implies, for all $\boldsymbol{m}\in\mathbb{N}_0^{p+1}$, $s_{\mu}(\boldsymbol{m})=s_{\nu}(\boldsymbol{m})$, hence $\mu=\nu$. 

To prove the second statement, assume there exists $P\in\R[Z_1,\dots,Z_p]$ nonzero of degree $k\in\N$, such that, $\forall \boldsymbol{x}\in\mathbb{S}_{\boldsymbol{X}}$, $P(\boldsymbol{x})=0$. Let $Q$ be the corresponding homogeneous polynomial in 
$\R[Z_1,\dots,Z_{p+1}]$ of degree $k$ such that, $\forall \boldsymbol{x}\in\mathbb{S}_{\boldsymbol{X}}$, $Q(1,\boldsymbol{x})=P(\boldsymbol{x})$. 
Let $\mathbb{P}_{\boldsymbol{\Gamma}}^*$ of support in $[-1,1]^{p+1}$ and density with respect to the Lebesgue measure bounded from below by $c>0$, and the infinitely differentiable  function $q(\boldsymbol{\gamma}) = \prod_{j=1}^{p+1} \boldsymbol{\gamma}_j \exp(-1/(\boldsymbol{\gamma}_j^2-1)) \indic\{\boldsymbol{\gamma}_j\in [-1,1]\}$. Then,
$\mathbb{P}_{\boldsymbol{\Gamma}} = \mathbb{P}_{\boldsymbol{\Gamma}}^* + c h /\|h\|_{\infty} d\boldsymbol{\gamma}$, where $h = Q(\partial_1, \dots, \partial_{p+1}) q $ is a probability. By properties of the Fourier transform of measures in $\mathfrak{M}^*(\R^{p+1})$ and homogeneity, we conclude from
$$
\mathcal{F}\left[\mathbb{P}_{\boldsymbol{\Gamma}} \right](t,t\boldsymbol{x}) = \mathcal{F}\left[\mathbb{P}_{\boldsymbol{\Gamma}}^* \right](t,t\boldsymbol{x})  + t^k Q(1,\boldsymbol{x})\mathcal{F}\left[h\right](t,t\boldsymbol{x}).
$$
\end{proof}

\begin{proof}[\bf Proof of Proposition \ref{prop:basis}.]
Start by proving (P6.\ref{P11}). Take $J\subseteq \Z\setminus\{0\}$ finite , $J_s\subseteq J$ the positive indices $j$ in $J$ such that $-j\in J$.\\
Let $(b_j)_{j\in J} \in \mathbb{N}_0^{|J|}$ such that $ \sum_{j\in J} b_j \lambda_j = 0$. We have
\begin{equation*}
\sum_{j\in J} b_j \lambda_j = 
\sum_{j\in J_s}  \frac{b_j+b_{-j}}{4r^{j}} + \sum_{j\in J\setminus (J_s\cup -J_s)}\frac{b_j}{4r^{|j|}} + \sum_{j\in J} j b_j\label{eqeric},
\end{equation*}
hence, $P(r)=0$, where $j_0 = \max_{j\in J} \abs{j}$ and
$$P(X)=\sum_{j\in J_s}  (b_j+b_{-j})X^{j_0-j}+
\sum_{j\in J\setminus J_s}b_jX^{j_0-|j|} + 4X^{j_0}\sum_{j\in J} j b_j \in\mathbb{Z}[X].$$
Because $r$ is transcendental over $\mathbb{Z}$ and $0\notin J$, for all $j\in J\setminus (J_s\cup -J_s)$, 
$b_j =0$, and, for all $j\in J_s$, $b_j+b_{-j}=0$ thus $b_j=b_{-j}=0$ because $b_j,b_{-j}\in\N_0$. Hence, for all $j\in J$, $b_j=0$.\\
Now, take $P\in \C(X)[X_1,\dots,X_p]$, where $p=|J|$, which is zero when evaluated at $(f_j)_{j\in J}$. Hence, we have 
$\sum_k c_k(z)\exp(z\sum_{j\in J}b_{k,j}\lambda_j)=0$,
where the sum over $k$ is finite, all $b_{k,j}$ belong to $\N_0$, and $c_k$ are rational functions. All exponentials in this sum are distinct by the above computations and we conclude that all $c_k$s are zero by taking limits in $\C$. This yields the result.\\
(P6.\ref{P12}) follows from Kadec's $1/4$-Theorem (see Theorem 14 page 42 in \cite{Young}) because $ \sup_{j\in\Z}|\lambda_j - j|<1/4$.
\end{proof}

\begin{proof}[\bf Proof of Theorem \ref{tL} under $\mathcal{R}_{L,j}$ for $j\in\{4,5\}$.] Denote  by  $ <,>_{\ell_2(\N)}$ the inner product in $\ell_2(\N)$. Start by considering $\mathcal{R}_{L,4}$. 
By a similar argument as the one involving \eqref{eq:repar},  we can consider without loss of generality $\{\underline{x}_m\} = \{0_m\}$ in $\mathcal{R}_{L,4}$ \eqref{RL3i}. Take $S = \ell_2(\N)$ and $\mathfrak{F}(S)$ the vector space spanned by $\mathcal{F}[\mathcal{P}]$.  An element of $\mathfrak{F}(S)$ is of the form $\mathcal{F}[\mu-\nu]$, where $\mu,\nu\in\mathcal{P}$.
	Let $n\in \N$. 
	We have, for all $\{z_m\} \in G_{n+1}$ and $t\in \R$, 
	$	\int_{\ell_2(\N)} e^{it <\{\gamma_m\},\{z_m\}>_{\ell_2(\N)}} d(\mu-\nu)(\{\gamma_m\})  
	= \int_{ G_{n+1}} e^{it <\{\gamma_m\},\{z_m\}>_{\ell_2(\N)}} d \Pi_{ G_{n+1}*}(\mu-\nu)(\{\gamma_m\}) $. 	We can now proceed like in the proof of Theorem \ref{tL} under $\mathcal{R}_{L,3}$ and obtain $\Pi_{G_{n+1}*}\mu=\Pi_{G_{n+1}*}\nu$ hence
	$ \mathcal{F}\left[ \Pi_{G_{n+1}*}\mu- \Pi_{G_{n+1}*}\nu\right]=0  $. Using that $\cup_{n\in \N} G_{n+1}$ is dense in $\ell_2(\N)$ and the continuity of the Fourier transform we conclude $ \mu = \nu$.\\
Identification under $\mathcal{R}_{L,5}$ is corollary of Theorem 4.1 in \cite{cuesta2007sharp}. 
\end{proof}

\begin{proof}[\bf Proof of Theorem \ref{tL} under $\mathcal{R}_{L,7}$.] The proof relies on a weaker form of diagram \eqref{ediag}. 
Take $\mathbb{P}_{\boldsymbol{\Gamma}}\in\mathcal{P}$ and $U_{l}  = \prod_{k=1}^l V_{k}$ for all $l\in[p]$. For $\boldsymbol{\beta} \in \C^p$ and $j\in[p]$, denote by  $\boldsymbol{\beta}_{[j]}=(\boldsymbol{\beta}_1,\dots,\boldsymbol{\beta}_j)^{\top}$. The function 
$\mathcal{F}[\mathbb{P}_{\boldsymbol{\Gamma}}-\mathbb{P}_{\boldsymbol{\Gamma}}^*]$ is 0 on $U$, so for 
 $t\ne0$ and $\boldsymbol{u}\in U_{p-1}$, $g_{p,t,\boldsymbol{u}}: z\in\R\mapsto  \mathcal{F}[\mathbb{P}_{\boldsymbol{\Gamma}}-\mathbb{P}_{\boldsymbol{\Gamma}}^*](t,t\mt{u},z)=\mathcal{F}[\mu_{p,t,\boldsymbol{u}}(\cdot)](z)$ is 0 on $tV_{p}$, where 
	$$ \mu_{p,t,\boldsymbol{u}}(\cdot) =\int_{\R^{p}}e^{it\left(a+\boldsymbol{b}_{[p-1]}^{\top}\boldsymbol{u}\right)}d(\mathbb{P}_{\boldsymbol{\Gamma}}-\mathbb{P}_{\boldsymbol{\Gamma}}^*)\left(a,\boldsymbol{b}_{[p-1]},\cdot\right).$$  
For $h$ nonnegative and measurable,  we have 
	$\int_{\R}h\left(z\right)d|\mu_{p,t,\boldsymbol{u}}|(z)\le\mathbb{E}_{\mathbb{P}_{\boldsymbol{\Gamma}}}\left[h\left(\boldsymbol{\beta}_p\right)\right]+\mathbb{E}_{\mathbb{P}_{\boldsymbol{\Gamma}}^*}\left[h\left(\boldsymbol{\beta}_p\right)\right]\le 2 M_k(h)$ (see Theorem 6.2 in \cite{Rudin2}), so $g_{p,t,\boldsymbol{u}}\in C_p$. Because $tV_{p}$ is a set of uniqueness of $C_p$, $g_{p,t,\boldsymbol{u}}$ is 0 on $\R$.\\ 	
	Take now $(\boldsymbol{u},x)\in U_{p-2}\times \R$ and $g_{p-1,t,\boldsymbol{u},x}: z\in\R\mapsto  \mathcal{F}[\mathbb{P}_{\boldsymbol{\Gamma}}-\mathbb{P}_{\boldsymbol{\Gamma}}^*](t,t\mt{u},z,x)=\mathcal{F}[\mu_{p-1,t,\boldsymbol{u},x}(\cdot)](z)$ is 0 on $tV_{p-1}$, where 
	$$\mu_{p-1,t,\boldsymbol{u},x}(\cdot)=\int_{\R^{p}}e^{it\left(a+\boldsymbol{b}_{[p-2]}^{\top}\boldsymbol{u}_{[p-2]}+\boldsymbol{b}_{p}x\right)}d\mu\left(a,\boldsymbol{b}_{[p-2]},\cdot,\boldsymbol{b}_{p}\right)$$
	Because $g_{p-1,t,\boldsymbol{u},x}\in C_{p-1}$ and $tV_{p-1}$ is set of uniqueness, $g_{p-1,t,\boldsymbol{u},x}$ is identically 0. We conclude by iterating this argument. 
	\end{proof}

\begin{proof}[\bf Proof of Theorem \ref{theo:binary} under $\mathcal{R}_{BC,1}$.]
Take $0<\epsilon\le1/2$. Assume that $\mathbb{P}_{\mt{\Gamma}}$ and $\mathbb{P}_{\mt{\Gamma}}^*$ both give rise to the same collection $ (\E\left[\indic\{\mt{\Gamma}^{\top}\mt{s}\geq 0\} \middle| \mt{S} = \mt{s} \right])_{\mt{s}\in\mathbb{S}_{\mt{S}}} $.
 Recall that, from \cite{GK},
\begin{equation}\label{eodd}
\mathcal{T}f(\boldsymbol{s})=\mathcal{T}f^-(\boldsymbol{s})=\left(\mathcal{T}f\right)^-(\boldsymbol{s}).
\end{equation}
One has in $S'\left(\mathbb{S}^p\right)$,  for all $f\in L^q\left(\mathbb{S}^p\right)$,
\begin{equation}\label{eegb}
\mathcal{T} f=\sum_{m\in\N_0}\lambda_{2m+1,p}Q_{2m+1,p}f,
\end{equation}
where $\lambda(2m+1,p)=(-1)^m2\pi^{p/2}1\cdot3\cdots(2m-1)/(\Gamma(p/2)p(p+2)\cdots(p+2m))$ for all $m\in\N_0$. 
Hence, we have, for all $\boldsymbol{x}\in\mathbb{S}_{\boldsymbol{X}}$,
$$\mathbb{E}[Y|\boldsymbol{X}=\boldsymbol{x}]=\sum_{m\in\N_0}\lambda_{2m+1,p}Q_{2m+1,p}f_{\boldsymbol{\Gamma}}\left(\frac{(1,\boldsymbol{x})}{\sqrt{1+|\boldsymbol{x}|_2^2}}\right)+\frac12$$
and similarly for $f_{\boldsymbol{\Gamma}}^*$. Denote the dot product and Lie norm by $\boldsymbol{z}^2=\sum_{k=1}^p\boldsymbol{z}_k^2$ and $L(\boldsymbol{z})=\sqrt{|\boldsymbol{z}|_2^2+\sqrt{|\boldsymbol{z}|_2^4-|\boldsymbol{z}^2|^2}}$ respectively. Denote also by
$G(\boldsymbol{z})=\sum_{m\in\N_0}\lambda_{2m+1,p}Q_{2m+1,p}f_{\boldsymbol{\Gamma}}(\boldsymbol{z}) +1/2$ for all $\boldsymbol{z}\in\C^{p+1}$ (resp. $G^*$ for $f_{\boldsymbol{\Gamma}}^*$)  and  
$F(\boldsymbol{z})=(f_1(\boldsymbol{z}),\dots,f_{p+1}(\boldsymbol{z}))^{\top}$ for all $\boldsymbol{z}\in\C^{p}$, where:
$$
f_1(\boldsymbol{z})=\frac{1}{\sqrt{1+\boldsymbol{z}^2}},\ 
f_2(\boldsymbol{z})=\frac{\boldsymbol{z}_1}{\sqrt{1+\boldsymbol{z}^2}},\ \dots\ ,\ 
f_{p+1}(\boldsymbol{z})=\frac{\boldsymbol{z}_p}{\sqrt{1+\boldsymbol{z}^2}}.
$$
We have $\mathbb{E}[Y|\boldsymbol{X}=\boldsymbol{x}]=G\circ F(\boldsymbol{x})=G^*\circ F(\boldsymbol{x})$ and we  prove  that $G\circ F, G^*\circ F\in\mathcal{A}(\R^p+i\epsilon\mathbb{B}_{\R}^p)$, hence $(G-G^*)\circ F\in\mathcal{A}(\R^p+i\epsilon\mathbb{B}_{\R}^p)$. 
For this, we check the conditions of Theorem 1.2.3 in \cite{RudinCn}. First, 
we have $F\in \mathcal{A}(\R^p+i\epsilon\mathbb{B}_{\R}^p)$. Indeed, for all $\boldsymbol{z}\in\R^p+i\epsilon\mathbb{B}_{\R}^p$, 
$1+\boldsymbol{z}^2\in(\C\setminus(-\infty,0])^{p+1}$. Now, for all $\boldsymbol{z}\in\R^p+i\epsilon\mathbb{B}_{\R}^p$, we have
\begin{equation*}
L(F(\boldsymbol{z}))^2=\frac{L(1,\boldsymbol{z})^2}{\left|\sqrt{1+\boldsymbol{z}^2}\right|^2}\le\frac{1+\left(|\mathrm{Re}(\boldsymbol{z})|_2+|\mathrm{Im}(\boldsymbol{z})|_2\right)^2}{1+|\mathrm{Re}(\boldsymbol{z})|_2^2-|\mathrm{Im}(\boldsymbol{z})|_2^2}\le1+2\epsilon.
\end{equation*}
Moreover, using Lemma 3.23 in \cite{Morimoto} for the first display and using the Young inequalities (see \cite{GK}) for the third display, we have, for all $\boldsymbol{z}\in \C^{p+1}$ such that $L(\boldsymbol{z})\le\sqrt{1+2\epsilon}$,
\begin{align*}\left|\lambda_{2m+1,p}Q_{2m+1,p}f_{\boldsymbol{\Gamma}}(\boldsymbol{z})\right|&\le L(\boldsymbol{z})^{2m+1}\left\|\lambda_{2m+1,p}Q_{2m+1,p}f_{\boldsymbol{\Gamma}}\right\|_{L^{\infty}\left(\mathbb{S}^p\right)}\\
&\le L(\boldsymbol{z})^{2m+1}\left\|\mathcal{T}\left[Q_{2m+1,p}f_{\boldsymbol{\Gamma}}\right]\right\|_{L^{\infty}\left(\mathbb{S}^p\right)}\\
&\le L(\boldsymbol{z})^{2m+1}\left\|Q_{2m+1,p}f_{\boldsymbol{\Gamma}}\right\|_{L^{1}\left(\mathbb{S}^p\right)}\\
& \le\left(1+2\epsilon\right)^{m+1/2}\left\|Q_{2m+1,p}f_{\boldsymbol{\Gamma}}\right\|_{L^{1}\left(\mathbb{S}^p\right)}.
\end{align*}
Now, because $\overline{\lim}_{m\to\infty}\left\|Q_{2m+1,p}f_{\boldsymbol{\Gamma}}\right\|_{L^{1}\left(\mathbb{S}^p\right)}^{1/m}\le q/(1+2\epsilon)$ with $q<1$, there exists $m_0\in\N$ such that, for all $m\ge m_0$,
$\left(1+2\epsilon\right)^{m+1/2}\left\|Q_{2m+1,p}f_{\boldsymbol{\Gamma}}\right\|_{L^{1}\left(\mathbb{S}^p\right)}\le\sqrt{1+2\epsilon}$.
Hence, there exists $C_{\epsilon}<\infty$ such that, for all $m\in\N_0$, 
$\left(1+2\epsilon\right)^{m+1/2}\left\|
Q_{2m+1,p}f_{\boldsymbol{\Gamma}}\right\|_{L^{1}\left(\mathbb{S}^p\right)}\le C_{\epsilon}$.
As a result, for all $\boldsymbol{z}\in \C^{p+1}$ such that $L(\boldsymbol{z})<\sqrt{1+2\epsilon}$, we have $\sup_{m\in\N_0}\left|\lambda_{2m+1,p}Q_{2m+1,p}f_{\boldsymbol{\Gamma}}(\boldsymbol{z})\right|<\infty$. Using the fact that $\overline{\boldsymbol{z}}\mapsto Q_{2m+1,p}f_{\boldsymbol{\Gamma}}(\overline{\boldsymbol{z}})$ are homogeneous harmonic polynomials and  Theorem 1.5.6 in \cite{RudinCn}
yields $G\circ F\in\mathcal{A}(\R^p+i\epsilon\mathbb{B}_{\R}^p)$. Moreover, for all $\boldsymbol{z}\in \R^p+i\epsilon\mathbb{B}_{\R}^p$, we have
\begin{align*}
\left|G\circ F(\boldsymbol{z})\right|&\le\sqrt{1+2\epsilon}\sum_{m\in\N_0}\left(\left(1+2\epsilon\right)\left\|
Q_{2m+1,p}f_{\boldsymbol{\Gamma}}\right\|_{L^{1}\left(\mathbb{S}^p\right)}^{1/m}\right)^m + \df{1}{2},
\end{align*}
and the upper bound is a convergent series  using $\mathcal{R}_{BC,1} \ \text{(ii)}$. Hence, $G\circ F\in\mathcal{A}( \R^p+i\epsilon\mathbb{B}_{\R}^p)$ is bounded and similarly for $G^*\circ F$, which yields the result.
\end{proof}

\begin{proof}[\bf Proof of Theorem \ref{theo:binary} under $\mathcal{R}_{BC,2}$.]
Assume that $\mathbb{P}_{\mt{\Gamma}}$ and $\mathbb{P}_{\mt{\Gamma}}^*$ both give rise to the same collection $ (\E\left[\mt{\Gamma}^{\top}\mt{s}\geq 0\} \middle| \mt{S} = \mt{s} \right])_{\mt{s}\in\mathbb{S}_{\mt{S}}} $.  
Due to \eqref{eegb},  
 $\mathcal{T}$ and $\Delta$ commute. For $r\ge0$, $q\in[1,\infty]$, and $f\in\ L^q\left(\mathbb{S}^p\right)$, $\Delta ^k f$ exists in $S'\left(\mathbb{S}^p\right)$ and is such that
$\Delta ^r f=\sum_{m=0}^{\infty}\zeta_{m,p}^{r}Q_{m,p}f$. 
When the sum converges in $L^q\left(\mathbb{S}^p\right)$,
$f$ is said to belong to $W_{q}^{2r}\left(\mathbb{S}^p\right)$. This Sobolev space is equipped with $\left\|f\right\|_{W_{q}^{2r}\left(\mathbb{S}^p\right)}=\left\|f\right\|_{L^q\left(\mathbb{S}^p\right)} +\left\|\Delta^rf\right\|_{L^q\left(\mathbb{S}^p\right)}$.  Using that $\Delta$ commutes with $\mathcal{T}$, \eqref{eodd}, the Young inequalities (see \cite{GK}), and $\mathcal{R}_{BC,2}$ \eqref{Rsphere5}, yield that, for all $m\in\N_0$, $\mathcal{T}f\in W_{\infty}^m\left(\mathbb{S}^p\right)$ and
\begin{equation}\label{emin}
\left\|\Delta^m \mathcal{T}f\right\|_{L^{\infty}\left(\mathbb{S}^p\right)}\le\left\|\Delta^mf^-\right\|_{L^1\left(\mathbb{S}^p\right)}.
\end{equation}
Denote, for $T\in\N$, by $\mathcal{K}_Tf(\cdot)=\int_{\mathbb{S}^p}k_T(\cdot,\boldsymbol{s})f(\boldsymbol{s})d\sigma(\boldsymbol{s})$, where
$k_T(\boldsymbol{x},\boldsymbol{y})=\sum_{l=0}^{T}\psi\left(l/T\right)q_{l,p}(\boldsymbol{x},\boldsymbol{y})$
with $\psi\in C^{\infty}([0,\infty))$, nonnegative, nonincreasing, such that $\psi(x)=1$ if $x\in[0,1]$, $0\le\psi(x)\le1$ if $x\in[1,2]$, and $\psi(x)=0$ if $x\ge2$, $\Delta$ also commutes with $\mathcal{K}_T$.
Now, using that $\Delta$ commutes with $\mathcal{T}$ and $\mathcal{K}_T$ (first and third displays) and the Young inequalities (first display), we obtain, for all $m\in\N_0$ and $k\in\N$,
\begin{align}
\left\|\Delta^m\mathcal{T}f-\mathcal{T}\mathcal{K}_T\Delta^m f\right\|_{\infty}
&\le\left\|\Delta^mf^--\Delta^m\mathcal{K}_Tf^-\right\|_{L^{1}\left(\mathbb{S}^p\right)}\notag\\
&=\left\|\Delta^mf^--\mathcal{K}_T\Delta^mf^-\right\|_{L^{1}\left(\mathbb{S}^p\right)}\notag\\
&\le CT^{-2k}\left\|\Delta^mf^-\right\|_{W_1^{2k}\left(\mathbb{S}^p\right)}\notag\\
&\le CT^{-2k}\left(\left\|f^-\right\|_{W_1^{2m}\left(\mathbb{S}^p\right)}+\left\|f^-\right\|_{W_1^{2m+2k}\left(\mathbb{S}^p\right)}\right)\label{ezero},
\end{align}
where the third display follows from Proposition A.2 in \cite{GK}. Using 
$$\mathcal{T}\mathcal{K}_T\Delta^m f=\sum_{l=0}^{\lfloor(T-1)/2\rfloor}\psi\left(\frac{2l+1}{T}\right)\lambda_{2l+1,p}Q_{2l+1,p}\Delta^m f,$$ that $q_{2l+1,p}(x,y)$ is a polynomial in $\boldsymbol{x}^{\top}\boldsymbol{y}$, and that $\Delta^m f\in L^1\left(\mathbb{S}^p\right)$, we obtain that $\mathcal{T}\mathcal{K}_T\Delta^m f$ is odd and continuous and conclude from \eqref{ezero} and  $\mathcal{R}_{BC,2}$ \eqref{Rsphere5} that $\mathcal{T}f\in C^{\infty}_{{\rm odd}}(\mathbb{S}^p)$.  $\mathcal{R}_{BC,2}$ \eqref{Rsphere5}  and \eqref{emin} yield 
\begin{align*}
\left\|\Delta^m \mathcal{T}\left(f_{\boldsymbol{\Gamma}}-f_{\boldsymbol{\Gamma}}^*\right)\right\|_{L^{\infty}\left(\mathbb{S}^p\right)}\le \widetilde{M}_m,
\end{align*}
where $ \widetilde{M}_m=\max(c_{f_{\boldsymbol{\Gamma}}}+c_{f_{\boldsymbol{\Gamma}}^*},1)\max(b_{f_{\boldsymbol{\Gamma}}},b_{f_{\boldsymbol{\Gamma}}^*})^m M_m$ is such that
$$\sum_{m\in\N}\widetilde{M}_m^{1/m}\ge  \max\left(c_{f_{\boldsymbol{\Gamma}}}+c_{f_{\boldsymbol{\Gamma}}^*},1\right)\max\left(b_{f_{\boldsymbol{\Gamma}}},b_{f_{\boldsymbol{\Gamma}}^*}\right)\sum_{m\in\N}M_m^{1/m}=\infty$$
because $\{M_m\}$ satisfies \eqref{ideni}. Also, \eqref{eeg}, the fact that $\mathcal{T}\left(f_{\boldsymbol{\Gamma}}-f_{\boldsymbol{\Gamma}}^*\right)\in C^{\infty}_{{\rm odd}}\left(\mathbb{S}^p\right)$, and  $\mathcal{R}_{BC,2}$ \eqref{Rsphere1} imply that, for all $\boldsymbol{u}\in U$ and $m\in\N_0$, $\Delta^m\mathcal{T}\left(f_{\boldsymbol{\Gamma}}-f_{\boldsymbol{\Gamma}}^*\right)(\boldsymbol{u})=0$.
From the proof of Theorem 10 in \cite{Bochner} (paragraph after Lemma 5) we obtain, for all $\boldsymbol{u}\in U$ and $l\in\N_0$,
$ Q_{2l+1,p}\mathcal{T}\left(f_{\boldsymbol{\Gamma}}-f_{\boldsymbol{\Gamma}}^*\right)(\boldsymbol{u})=0$.
Hence, by  $\mathcal{R}_{BC,2}$ \eqref{Rsphere1} and the fact that $Q_{2l+1,p}\mathcal{T}\left(f_{\boldsymbol{\Gamma}}-f_{\boldsymbol{\Gamma}}^*\right)$ is the restriction to $\mathbb{S}^p$ of a harmonic homogeneous polynomial of degree $2l+1$, we obtain, for all  $\boldsymbol{s}\in \mathbb{S}^p$ and  $l\in\N_0$,
$Q_{2l+1,p}\mathcal{T}\left(f_{\boldsymbol{\Gamma}}-f_{\boldsymbol{\Gamma}}^*\right)(\boldsymbol{s})=0$.
Thus $\left(f_{\boldsymbol{\Gamma}}-f_{\boldsymbol{\Gamma}}^*\right)^-=0$ and the odd parts of $f_{\boldsymbol{\Gamma}}$ and $f_{\boldsymbol{\Gamma}}^*$ coincide in $L^{1}\left(\mathbb{S}^p\right)$. \\
By the argument in \cite {GLP}, using that for $f=f_{\boldsymbol{\Gamma}}$ or $f=f_{\boldsymbol{\Gamma}}^*$, for a.e. $\boldsymbol{u}\in\mathbb{S}^p$ $f(\boldsymbol{u})\ge0$, the definition of $f^{-}$, and 
$f_{\boldsymbol{\Gamma}}(\boldsymbol{u})f_{\boldsymbol{\Gamma}}(-\boldsymbol{u})=0$ for a.e. $\boldsymbol{u}\in\mathbb{S}^p$, we obtain that, for a.e. $\boldsymbol{u}\in\mathbb{S}^p$,
$f(\boldsymbol{u})=2 f^{-}(\boldsymbol{u})\indic\{f^{-}(\boldsymbol{u})>0\}$.
Thus, there is a one-to-one mapping between $f$ and $f^-$. This allows to conclude.
\end{proof}

\begin{proof}[\bf Proof of Theorem \ref{Kot}.]  
	Assume that $P=\mathbb{P}_{\boldsymbol{\Gamma}}$ and $P=\mathbb{P}_{\boldsymbol{\Gamma}}^*$ both give rise to the same pair	$\left(\mathbb{P}_{\boldsymbol{\Gamma}_1 + \boldsymbol{\Gamma}_3}
	, \mathbb{P}_{\boldsymbol{\Gamma}_2 + \boldsymbol{\Gamma}_3}\right)$.  
	For $t_0>0$ small enough, on $(-t_0,t_0)$, $\varphi_{\delta}$, $\varphi_{\delta}^*$, $\varphi_{\boldsymbol{\epsilon}_1}$, $\varphi_{\boldsymbol{\epsilon}_1}^*$, $\varphi_{\boldsymbol{\epsilon}_2}$, $\varphi_{\boldsymbol{\epsilon}_2}^*$ do not vanish, hence there exist nonvanishing continuous functions $p_{\delta}$, $p_{\boldsymbol{\epsilon}_1}$, and $p_{\boldsymbol{\epsilon}_2}$ such that  
	$\varphi_{\delta}=p_{\delta}\varphi_{\delta}^*$, $\varphi_{\boldsymbol{\epsilon}_1}=p_{\boldsymbol{\epsilon}_1}\varphi_{\boldsymbol{\epsilon}_1}^*$, $\varphi_{\boldsymbol{\epsilon}_2}=p_{\boldsymbol{\epsilon}_2}\varphi_{\boldsymbol{\epsilon}_2}^*$. By the restriction $\mathcal{R}_K$, we have $p_{\boldsymbol{\epsilon}_1}\in C^{\infty}(-t_0,t_0)$. 
	We have
	\begin{equation}\label{eK0}p_{\delta}(\boldsymbol{t}_1+\boldsymbol{t}_2)p_{\boldsymbol{\epsilon}_1}(\boldsymbol{t}_1)p_{\boldsymbol{\epsilon}_2}(\boldsymbol{t}_2)=1\quad \text{for all}\ -t_0<  \boldsymbol{t}_1,\boldsymbol{t}_2,\boldsymbol{t}_1+\boldsymbol{t}_2< t_0,
	\end{equation}
	hence
	\begin{equation}\label{eK1}
	p_{\delta}(t)p_{\boldsymbol{\epsilon}_1}(t)=1\quad \text{and}\quad p_{\delta}(t)p_{\boldsymbol{\epsilon}_2}(t)=1\quad \text{for  all}\ -t_0< t< t_0.
	\end{equation}
	Injecting \eqref{eK1} into \eqref{eK0} we obtain, for all $-t_0< \boldsymbol{t}_1,\boldsymbol{t}_2,\boldsymbol{t}_1+\boldsymbol{t}_2< t_0$,
	$p_{\delta}(\boldsymbol{t}_1+\boldsymbol{t}_2)=p_{\delta}(\boldsymbol{t}_1)p_{\delta}(\boldsymbol{t}_2) $, 
	which, using again \eqref{eK1}, yields
	$p_{\boldsymbol{\epsilon}_1}(\boldsymbol{t}_1+\boldsymbol{t}_2)=p_{\boldsymbol{\epsilon}_1}(\boldsymbol{t}_1)p_{\boldsymbol{\epsilon}_1}(\boldsymbol{t}_2)$. 	Hence, we have 
	$p_{\boldsymbol{\epsilon}_1}'(\boldsymbol{t}_1+\boldsymbol{t}_2)=p_{\boldsymbol{\epsilon}_1}'(\boldsymbol{t}_1)p_{\boldsymbol{\epsilon}_1}(\boldsymbol{t}_2)$,
	which at $\boldsymbol{t}_1=0$ and $\boldsymbol{t}_2=t$ yields
	\begin{equation}\label{edoexp}
	p_{\boldsymbol{\epsilon}_1}'(t)=p_{\boldsymbol{\epsilon}_1}'(0)p_{\boldsymbol{\epsilon}_1}(t)\quad \text{for all} \ -t_0< t< t_0,
	\end{equation}
	where
	$p_{\boldsymbol{\epsilon}_1}'(0)=\varphi_{\boldsymbol{\epsilon}_1}'(0)-\left(\varphi_{\boldsymbol{\epsilon}_1}^*\right)'(0)=i(\E_{\mathbb{P}_{\boldsymbol{\epsilon}_1}}[\boldsymbol{\epsilon}_1]-\E_{\mathbb{P}_{\boldsymbol{\epsilon}_1}^*}[\boldsymbol{\epsilon}_1]) \in i\R$
	which we denote by 
	$p_{\boldsymbol{\epsilon}_1}'(0):=ib$.
	Thus, we obtain,  for all $t_0< t< t_0$, $p_{\boldsymbol{\epsilon}_1}(t)=\exp(ibt)$. This yields that  $t\mapsto \varphi_{\boldsymbol{\epsilon}_1}(t)-\exp(ibt)\varphi_{\boldsymbol{\epsilon}_1}^*(t)$ is 0 on $(-t_0,t_0)$. Thus, using $\mathcal{R}_K$, we have $\varphi^{(1)}_{\boldsymbol{\epsilon}_1}(0)=ib +(\varphi^*_{\boldsymbol{\epsilon}_1})^{(1)}(0)$ then using that $\mathbb{P}_{\boldsymbol{\epsilon}_1}$ and $\mathbb{P}_{\boldsymbol{\epsilon}_1}^*$ are both mean 0 yields $b=0$. Thus, $\mathbb{P}_{\boldsymbol{\epsilon}_1}$ and $\mathbb{P}_{\boldsymbol{\epsilon}_1}^*$ have the same moments so $\mathbb{P}_{\boldsymbol{\epsilon}_1} =\mathbb{P}_{\boldsymbol{\epsilon}_1}^*$ by $\mathcal{R}_K$. 

	Now, for all $t\in\R$, $\varphi_{\boldsymbol{Y}_1}(t)=\varphi_{\delta}(t)\varphi_{\boldsymbol{\epsilon}_1}^*(t)=\varphi_{\delta}^*(t)\varphi_{\boldsymbol{\epsilon}_1}^*(t)$,  hence, because the zeros of $\varphi_{\boldsymbol{\epsilon}_1}^*$ are isolated (see Lemma 4.8 in \cite{Belisle}) and $\varphi_{\delta}$ and $\varphi_{\delta}^*$ are continuous, we obtain $\varphi_{\delta}=\varphi_{\delta}^*$.\\
	Similarly, because, for all $t\in\R$, $\varphi_{\boldsymbol{Y}_2-\boldsymbol{Y}_1}(t)=\varphi_{\boldsymbol{\epsilon}_2}^*(t)\varphi_{\boldsymbol{\epsilon}_1}^*(-t)=\varphi_{\boldsymbol{\epsilon}_2}(t)\varphi_{\boldsymbol{\epsilon}_1}^*(-t)$, the zeros of $\varphi_{\boldsymbol{\epsilon}_1}^*$ are isolated and $\varphi_{\boldsymbol{\epsilon}_2}$ and $\varphi_{\boldsymbol{\epsilon}_2}^*$ are continuous, we obtain $\varphi_{\boldsymbol{\epsilon}_2}=\varphi_{\boldsymbol{\epsilon}_2}^*$,
	hence the result.
	\end{proof}

\begin{proof}[\bf Proof of Theorem \ref{theo:LP2T}.]
	We use $\mathcal{X}_0= \{ \mt{x}\in\R^T:\ \prod_{j=1}^T\mt{x}_j\ne0,\   \prod_{m\neq j} (\mt{x}_m - \mt{x}_j)\ne0\}$ and, for all $(\mt{v},\mt{x})\in \R^T\times\mathcal{X}_0$,
	\begin{equation}\label{eq:inv_V}
	\Theta:\ (\mt{v},\mt{x}) \mapsto \sum_{k=1}^T  \mt{v}_k \left(\sum_{j=1}^T b_{jk}(\mt{x}) \boldsymbol{x}_j^{T}\right),
	\end{equation}
	\begin{align*}\left\{ \begin{array}{rll}
	b_{jk}(\mt{x}) = & \displaystyle\dfrac{(-1)^{k+1}}{\displaystyle\prod_{\underset{m\neq j}{1\leq m \leq T}} (\mt{x}_m - \mt{x}_j)}  \sum_{\underset{i_1,\dots,i_{T-k}\neq j}{1\leq i_1 < \dots < i_{T-k}\leq T}} \mt{x}_{i_1}\mt{x}_{i_2} \dots \mt{x}_{i_{T-k}} & \text{for all} \ k \neq T  \\ 
	b_{jk}(\mt{x}) = & \dfrac{(-1)^{T+1}}{\displaystyle\prod_{\underset{m\neq j}{1\leq m \leq T}} (\mt{x}_m - \mt{x}_j)} & \text{for } \ k = T .
	\end{array} \right.
	\end{align*}
	Assume that $\mathbb{P}_{\boldsymbol{\Gamma}}$ and $\mathbb{P}_{\boldsymbol{\Gamma}}^*$ both give rise to the same collections	$\left(\mathbb{P}_{\boldsymbol{\Gamma}_1 + \sum_{j=1}^T\mt{\Gamma}_{j+1}\mt{x}^j_t + \epsilon_t}\right)_{\mt{x}_t \in  \mathbb{S}_{\mt{X}_t}, t=1,\dots, T}$.  
	Using  ${\mathcal{R}_{LP,0}}$ \eqref{R0LPii}, there exists $r>0$ such that $r \mt{1}_c \in \mathcal{X}_1 $ and conditioning on $\boldsymbol{X} = r \mt{1}_c$ and using $\delta =\alpha+\sum_{k=1}^T\boldsymbol{\beta}_k r^k$ yield, for all $\mt{t} \in \R^T$, 
	$\E\left[ e^{i\mt{t}^{\top}\mt{Y}} \middle| \boldsymbol{X} = r \mt{1}_c \right]= \varphi_{\delta}\left(\sum_{j=1}^{T}\mt{t}_j\right)\prod_{j=1}^T\varphi_{ \boldsymbol{\epsilon}_j}(\mt{t}_j)$,
	hence 
	$$ \varphi_{\delta}\left(\sum_{j=1}^{T}\mt{t}_j\right)\prod_{j=1}^T\varphi_{ \boldsymbol{\epsilon}_j}(\mt{t}_j) = \varphi_{\delta}^*\left(\sum_{j=1}^{T}\mt{t}_j\right)\prod_{j=1}^T\varphi^*_{ \boldsymbol{\epsilon}_j}(\mt{t}_j).$$
	Then, following the same steps as in the proof of Theorem \ref{Kot} yields $\mathbb{P}_{\boldsymbol{\epsilon}_1}=\mathbb{P}_{\boldsymbol{\epsilon}_1}^*$.  Now, for all $t\in\R$, 
	$$
	\varphi_{\delta}(t)\varphi_{\boldsymbol{\epsilon}_1}^*(t)=\varphi_{\delta}^*(t)\varphi_{\boldsymbol{\epsilon}_1}^*(t),$$  hence, because the zeros of $\varphi_{\boldsymbol{\epsilon}_1}^*$ are isolated (see Lemma 4.8 in \cite{Belisle}) and $\varphi_{\delta}$ and $\varphi_{\delta}^*$ are continuous, we obtain $\varphi_{\delta}=\varphi_{\delta}^*$. Similarly, because, for all $t\in\R$ and $j=2,\dots, T$, 
	$\E\left[e^{it(\boldsymbol{Y}_j-\boldsymbol{Y}_{1})}| (\mt{X}_{1}, \mt{X}_{j} )= (r,r)\right]=\varphi_{\boldsymbol{\epsilon}_j}^*(t)\varphi_{\boldsymbol{\epsilon}_1}^*(-t)$,
	we have $$\varphi_{\boldsymbol{\epsilon}_j}^*(t)\varphi_{\boldsymbol{\epsilon}_1}^*(-t)=\varphi_{\boldsymbol{\epsilon}_j}(t)\varphi_{\boldsymbol{\epsilon}_1}^*(-t).$$ Thus, using that the zeros of $\varphi_{\boldsymbol{\epsilon}_1}^*$ are isolated and $\varphi_{\boldsymbol{\epsilon}_j}$ and $\varphi_{\boldsymbol{\epsilon}_j}^*$ are continuous for all $j=2,\dots, T$, we obtain $\varphi_{\boldsymbol{\epsilon}_j}=\varphi_{\boldsymbol{\epsilon}_j}^*$ for all  $j=2,\dots, T$, hence, for all $\boldsymbol{t}\in\R^{T}$ and $\boldsymbol{x}\in \mathbb{S}_{\boldsymbol{X}}$,
	\begin{equation}\label{eq:V11}
	\varphi_{\alpha, \boldsymbol{\beta}}\left(\sum_{j=1}^{T}\mt{t}_j, \sum_{j=1}^{T} \mt{t}_j \mt{x}_j, \dots, \sum_{j=1}^{T} \mt{t}_j \mt{x}_j^{T} \right) =	\varphi_{\alpha, \boldsymbol{\beta}}^*\left(\sum_{j=1}^{T}\mt{t}_j, \sum_{j=1}^{T} \mt{t}_j \mt{x}_j, \dots, \sum_{j=1}^{T} \mt{t}_j \mt{x}_j^{T} \right).
	\end{equation} 
	Then, for all $\boldsymbol{x}\in\mathcal{X}_0$, we use a change of variable that relates $\boldsymbol{t}\in\R^{T}$ to $\boldsymbol{v}\in\R^{T}$ such that $(\sum_{j=1}^{T}\mt{t}_j, \sum_{j=1}^{T} \mt{t}_j \mt{x}_j, \dots, \sum_{j=1}^{T-1} \mt{t}_j \mt{x}_j^{T-1} ) = \mt{v}$. This change of variable allows to choose the values of the $T-1$ first variables in \eqref{eq:V11} independently from each other. This can be written as $   \mt{t}  = \left(V^{-1}(\mt{x})\right)^{\top} \mt{v}$, where 
	$$ V(\mt{x})= \left( \begin{array}{ccccc}
	1 & \mt{x}_1 &  \mt{x}_1^2 & \dots  & \mt{x}_1^{T-1} \\ 
	\colon &  &  & \dots  & \colon \\ 
	1 & \mt{x}_T&  \mt{x}_T^2 & \dots  & \mt{x}_T^{T-1}
	\end{array}  \right)$$
	is the Vandermonde matrix.  We use $D(\mt{x})$, the diagonal matrix which entries are the coordinates of $\mt{x}$. It is a classical result that, for all $\boldsymbol{x}\in\mathcal{X}_0$, $D(\mt{x})$ and $V(\mt{x})$ are invertible hence $\widetilde{V}(\mt{x}) = (D(\mt{x}) V(\mt{x}))^{\top} $ is invertible. Then, for all $\mt{v}\in \R^T$ and $ \mt{x} \in \mathcal{X}_0$, we can express $\mt{t}$ as a function of $\mt{v}$ and $ \mt{x}$, and obtain for the last variable in 
	\begin{equation*}
	\varphi_{\alpha, \boldsymbol{\beta}}^*\left(\sum_{j=1}^{T}\mt{t}_j, \sum_{j=1}^{T} \mt{t}_j \mt{x}_j, \dots, \sum_{j=1}^{T} \mt{t}_j \mt{x}_j^{T} \right),
	\end{equation*}
	that
	\begin{align}
	\sum_{j=1}^{T} \mt{t}_j \mt{x}_j^T& = \mt{t}^{\top} \mt{x}^T = \left(\left(V(\mt{x})^{-1}\right)^{\top}\mt{v}\right)^{\top}\mt{x}^{T} =\left( D(\mt{x})\widetilde{V}(\mt{x})^{-1} \mt{v}\right)^{\top}\mt{x}^{T}=\Theta(\mt{v},\mt{x}).\label{eq:change11}
	\end{align}
	Using  \eqref{eq:V11}, this yields, for all $\mt{v} \in\R^{T}$ and $\boldsymbol{x}\in\mathcal{X}_0$, \begin{equation}\label{eq:chang}
	\varphi_{\alpha, \boldsymbol{\beta}}\left(\mt{v}, \Theta(\mt{v},\mt{x})  \right) = \varphi_{\alpha, \boldsymbol{\beta}}^*\left(\mt{v}, \Theta(\mt{v},\mt{x})  \right).
	\end{equation}  
	Then, using that the vector space spanned by $\mathcal{F}[\mathcal{P}_c(\Omega_T)]$ is a quasi-analytic class of functions on $\R$, we obtain that, for all $\mt{v} \in\R^{T}$, $ z \in \R\mapsto (\varphi_{\alpha,\boldsymbol{\beta}} - \varphi^*_{\alpha,\boldsymbol{\beta}})(\mt{v},z )$
	belongs to a quasi-analytic class. Finally, for all $\mt{v}\in\R^{T}\setminus \mathcal{V}$ where $\mathcal{V}$ is such that $\R^{T}\setminus \mathcal{V}$ is dense in $\R^{T}$,  $\mathcal{R}_{LP}$ \eqref{R1LPiiiT} and that $ \mt{x} \mapsto \Theta(\mt{v},\mt{x})$ is continuous on  $\mathbb{S}_{\boldsymbol{X}}  \cap  \mathcal{X}_0$ yield that $U_{T,\mt{v}}=\left\{ \Theta(\mt{v},\mt{u}),\ \forall \boldsymbol{u} \in \mathbb{S}_{\boldsymbol{X}}  \cap  \mathcal{X}_0\right\}$ contains a bounded sequence of distinct points. We obtain that, for all $\mt{v} \in\R^{T}\setminus \mathcal{V}$ and $ z \in \R$, 
	$ \varphi_{\alpha,\boldsymbol{\beta}}(\mt{v},z ) = \varphi^*_{\alpha,\boldsymbol{\beta}}(\mt{v},z )$
	hence $\mathbb{P}_{\alpha,\boldsymbol{\beta}} = \mathbb{P}_{\alpha,\boldsymbol{\beta}}^*$ by continuity for all $\mt{v}\in\R^{T}$. \end{proof}

\begin{proof}[\bf Proof of Theorem \ref{cor:SI}.]
Let $0< \epsilon < 1$. Assume that $\mathbb{P}_{\boldsymbol{\Gamma}}$ and $\mathbb{P}_{\boldsymbol{\Gamma}}^*$ both give rise to the same collections	$\left(\mathbb{P}_{f(\boldsymbol{\Gamma}^{\top}\mt{x}_t) + \eta_t}\right)_{\mt{x}_t \in  \mathbb{S}_{\mt{X}_t}, t=1,2}$.  
Denote by $\mathcal{X}_1=\{(\boldsymbol{x}_1,\boldsymbol{x}_2)\in\mathbb{S}_{\boldsymbol{X}_1,\boldsymbol{X}_2}: \boldsymbol{x}_1=\boldsymbol{x}_2\}$. Using $\mathcal{R}_{SI}$ \eqref{R0SIii},  there exists $r >0$ such that $r (\mt{1}_c, \mt{1}_c) \in \mathcal{X}_1$, hence, for all $\mt{t} \in \R^{2}$, using $\delta : = f(\boldsymbol{\Gamma}^{\top} (r\mt{1}_c))$, we have 
\begin{equation}\label{eq:SI2}
\varphi_{\delta}(\mt{t}_1 + \mt{t}_2)\varphi_{ \boldsymbol{\eta}_1}(\mt{t}_1)\varphi_{ \boldsymbol{\eta}_2}(\mt{t}_2)
= \varphi^*_{\delta}(\mt{t}_1 + \mt{t}_2)\varphi^*_{ \boldsymbol{\eta}_1}(\mt{t}_1)\varphi^*_{ \boldsymbol{\eta}_2}(\mt{t}_2).
\end{equation}  
Then, following the same steps as in the proof of Theorem \ref{Kot} yields $ \varphi_{\mt{\eta}_1} =  \varphi_{\mt{\eta}_1}^*$. 
Now, for all $t\in\R$,  $\varphi_{\boldsymbol{Y}_1|\mt{X}_1}(t|r\mt{1}_c)=\varphi_{\delta}(t)\varphi_{\boldsymbol{\eta}_1}^*(t)=\varphi_{\delta}^*(t)\varphi_{\boldsymbol{\eta}_1}^*(t)$,  hence, because the zeros of $\varphi_{\boldsymbol{\eta}_1}^*$ are isolated (see Lemma 4.8 in \cite{Belisle}) and $\varphi_{\delta}$ and $\varphi_{\delta}^*$ are continuous, we obtain $\varphi_{\delta}=\varphi_{\delta}^*$. Similarly, because, for all $t\in\R$, $\varphi_{\boldsymbol{Y}_2-\boldsymbol{Y}_{1}| \mt{X}_1, \mt{X}_2}(t |(r\mt{1}_c,r\mt{1}_c))=\varphi_{\boldsymbol{\eta}_2}^*(t)\varphi_{\boldsymbol{\eta}_1}^*(-t)=\varphi_{\boldsymbol{\eta}_2}(t)\varphi_{\boldsymbol{\eta}_1}^*(-t)$, the zeros of $\varphi_{\boldsymbol{\eta}_1}^*$ are isolated and $\varphi_{\boldsymbol{\eta}_2}$ and $\varphi_{\boldsymbol{\eta}_2}^*$ are continuous, we obtain $\varphi_{\boldsymbol{\eta}_2}=\varphi_{\boldsymbol{\eta}_2}^*$.\\
Thus, we obtain, for all $\mt{t} \in \R^{2}$ and $(\boldsymbol{x}_1,\boldsymbol{x}_2)\in \mathbb{S}_{\mt{X}_1,\mt{X}_2}$,
\begin{align}\label{eq:SI1}
\varphi_{f(\boldsymbol{\Gamma}^{\top}\boldsymbol{x}_1),f(\boldsymbol{\Gamma}^{\top}\boldsymbol{x}_2)}(\mt{t})=\varphi^*_{f(\boldsymbol{\Gamma}^{\top}\boldsymbol{x}_1),f(\boldsymbol{\Gamma}^{\top}\boldsymbol{x}_2)}(\mt{t}). 
\end{align}
This amount to study identification in
$\boldsymbol{Z}_t=f(\boldsymbol{\Gamma}^{\top}\boldsymbol{X}_t)$, $t=1,2$
and $\boldsymbol{Z}_2\ge \boldsymbol{Z}_1$ is equivalent to $\boldsymbol{\Gamma}^{\top}\left(\boldsymbol{X}_{2}-\boldsymbol{X}_{1}\right)\ge0$.
This is the binary choice model, hence  Theorem \ref{theo:binary} 
yields the result. 
\end{proof}

\bibliographystyle{abbrv}
\bibliography{bib_full}

\begin{thebibliography}{10}

\bibitem{Belisle}
C.~Belisle, J.-C. Mass{\'e}, and T.~Ransford.
\newblock When is a probability measure determined by infinitely many
  projections?
\newblock {\em Annals of Probability}, 25:767--786, 1997.

\bibitem{Beran1}
R.~Beran and P.~Hall.
\newblock Estimating coefficient distributions in random coefficient
  regressions.
\newblock {\em Annals of Statistics}, 20:1970--1984, 1992.

\bibitem{BM}
R.~Beran and W.~Millar.
\newblock Minimum distance estimation in random coefficient regression models.
\newblock {\em Annals of Statistics}, 22:1976--1992, 1994.

\bibitem{Berg}
C.~Berg.
\newblock {\em The multidimensional moment problem and semi-groups, In: Landau,
  H.J. (ed.) Moments in Mathematics}.
\newblock American Mathematical Society, 1992.

\bibitem{berry1995automobile}
S.~Berry, J.~Levinsohn, and A.~Pakes.
\newblock Automobile prices in market equilibrium.
\newblock {\em Econometrica}, pages 841--890, 1995.

\bibitem{berry2009nonparametric}
S.~T. Berry and P.~A. Haile.
\newblock Nonparametric identification of multinomial choice demand models with
  heterogeneous consumers.
\newblock Technical report, National Bureau of Economic Research, 2009.

\bibitem{BCR}
J.~Bochnak, M.~Coste, and M.-F. Roy.
\newblock {\em Real algebraic geometry}.
\newblock Springer, 1998.

\bibitem{Bochner}
S.~Bochner and A.~Taylor.
\newblock Some theorems on quasi-analyticity for functions of several
  variables.
\newblock {\em American Journal of Mathematics}, 61:303--329, 1939.

\bibitem{breunig2019}
C.~Breunig.
\newblock Varying random coefficient models.
\newblock {\em Journal of Econometrics}, 221(2):381--408, 2021.

\bibitem{deconv}
M.~Carrasco and J.-P. Florens.
\newblock A spectral method for deconvolving a density.
\newblock {\em Econometric Theory}, 27:546--581, 2011.

\bibitem{Cater}
F.~S. Cater.
\newblock Differentiable, nowhere analytic functions.
\newblock {\em American Mathematical Monthly}, 91:618--624, 1984.

\bibitem{CP}
I.~Chalendar and J.~R. Partington.
\newblock Multivariable approximate {C}arleman-type theorems for complex
  measures.
\newblock {\em Annals of Probability}, 35:384--396, 2007.

\bibitem{CH}
V.~Chernozhukov and C.~Hansen.
\newblock An iv model of quantile treatment effets.
\newblock {\em Econometrica}, 73:245--261, 2005.

\bibitem{cooprider2020panel}
J.~Cooprider, S.~Hoderlein, and A.~Meister.
\newblock A panel data estimator for the distribution and quantiles of marginal
  effects in nonlinear structural models with an application to the demand for
  junk food.
\newblock {\em Available at SSRN 3545485}, 2020.

\bibitem{cuesta2007sharp}
J.~A. Cuesta-Albertos, R.~Fraiman, and T.~Ransford.
\newblock A sharp form of the {C}ramer-{W}old theorem.
\newblock {\em Journal of Theoretical Probability}, 20:201--209, 2007.

\bibitem{DeJeu1}
M.~{D}e {J}eu.
\newblock Determinate multidimensional measures, the extended {C}arleman
  theorem and quasi-analytic weights.
\newblock {\em Annals of Probability}, 20:1205--1227, 2003.

\bibitem{DeJeu2}
M.~{D}e {J}eu.
\newblock Subspaces with equal closure.
\newblock {\em Constructive Approximation}, 20:93--157, 2004.

\bibitem{evdokimov2012some}
K.~Evdokimov and H.~White.
\newblock Some extensions of a lemma of {K}otlarski.
\newblock {\em Econometric Theory}, 28:925--932, 2012.

\bibitem{fox2017note}
J.~Fox.
\newblock A note on nonparametric identification of distributions of random
  coefficients in multinomial choice models.
\newblock {\em NBER working paper series}, 2017.

\bibitem{Note}
C.~Gaillac and E.~Gautier.
\newblock Estimates for the {SVD} of the truncated fourier transform on
  ${L}^2\left(\cosh(b\cdot)\right)$ and stable analytic continuation.
\newblock 2019.
\newblock Preprint \href{https://arxiv.org/abs/1905.11338}{arXiv:1905.11338v3}.

\bibitem{estimation}
C.~Gaillac and E.~Gautier.
\newblock Adaptive estimation in the linear random coefficients model when
  regressors have limited variation.
\newblock {\em Forthcoming, Bernoulli}, 2021.

\bibitem{Gautier}
E.~Gautier.
\newblock Relaxing monotonicity in endogenous selection models and application
  to surveys.
\newblock In {\em Advances in Contemporary Statistics and Econometrics}, 2021.

\bibitem{Gautier2}
E.~Gautier and S.~Hoderlein.
\newblock A triangular treatment effect model with random coefficients in the
  selection equation.
\newblock 2011,2015.
\newblock Preprint \href{http://arxiv.org/pdf/1109.0362v4}{arXiv:1109.0362v4}.

\bibitem{GK}
E.~Gautier and Y.~Kitamura.
\newblock Nonparametric estimation in random coefficients binary choice models.
\newblock {\em Econometrica}, 81:581--607, 2013.

\bibitem{GLP}
E.~Gautier and E.~Le~Pennec.
\newblock Adaptive estimation in the nonparametric random coefficients binary
  choice model by needlet thresholding.
\newblock {\em Electronic Journal of Statistics}, 12:277--320, 2018.

\bibitem{Hardy1}
G.~H. Hardy.
\newblock On {S}tieltjes' ``probl\`eme des moments".
\newblock {\em Messenger of Math.}, 46:175--182, 1917.

\bibitem{Hardy2}
G.~H. Hardy.
\newblock On {S}tieltjes' ``probl\`eme des moments".
\newblock {\em Messenger of Math.}, 47:84--91, 1918.

\bibitem{HS}
J.~Heckman and B.~Singer.
\newblock A method for minimizing the impact of distributional assumptions in
  econometric models for duration data.
\newblock {\em Econometrica}, 52:271--320, 1984.

\bibitem{HV}
J.~Heckman and E.~Vytlacil.
\newblock Structural equations, treatment effects, and econometric policy
  evaluation1.
\newblock {\em Econometrica}, 73:669--738, 2005.

\bibitem{Heyde}
C.~C. Heyde.
\newblock Some remarks on the moment problem.
\newblock {\em The Quarterly Journal of Mathematics}, 14:91--96, 1963.

\bibitem{Hirschman}
I.~I. Hirschman.
\newblock On the distribution of the zeros of functions belonging to certain
  quasi-analytic classes.
\newblock {\em American Journal of Mathematics}, 72:396--406, 1950.

\bibitem{hoderlein2014triangular}
S.~Hoderlein, H.~Holzmann, and A.~Meister.
\newblock The triangular model with random coefficients.
\newblock {\em Journal of Econometrics}, 201:144--169, 2017.

\bibitem{IT}
H.~Ichimura and T.~S. Thompson.
\newblock Maximum likelihood estimation of a binary choice model with random
  coefficients of unknown distribution.
\newblock {\em Journal of Econometrics}, 86:269--295, 1998.

\bibitem{IA}
G.~W. Imbens and J.~D. Angrist.
\newblock Identification and estimation of local average treatment effects.
\newblock {\em Econometrica}, 62:467--475, 1994.

\bibitem{Stoy}
C.~Kleiber and J.~Stoyanov.
\newblock Multivariate distributions and the moment problem.
\newblock {\em Journal of Multivariate Analysis}, 113:7--18, 2013.

\bibitem{kotlarski}
I.~Kotlarski.
\newblock On characterizing the {G}amma and the normal distribution.
\newblock {\em Pacific Journal of Mathematics}, 20:69--76, 1967.

\bibitem{Lin2}
G.~D. Lin.
\newblock On the moment problems.
\newblock {\em Statistics and Probability Letters}, 35:85--90, 1997.

\bibitem{Lin}
G.~D. Lin.
\newblock Recent developments on the moment problem.
\newblock 2017.
\newblock Preprint \href{http://arxiv.org/pdf/1703.01027}{arXiv:1703.01027v1}.

\bibitem{Mand}
S.~Mandelbrojt.
\newblock {\em S\'eries adh\'erentes, r\'egularisation des suites,
  applications}.
\newblock Gauthiers-Villars, 1952.

\bibitem{masten2017random}
M.~A. Masten.
\newblock Random coefficients on endogenous variables in simultaneous equations
  models.
\newblock {\em The Review of Economic Studies}, 85(2):1193--1250, 2017.

\bibitem{Meister07}
A.~Meister.
\newblock Deconvolving compactly supported densities.
\newblock {\em Mathematical Methods of Statistics}, 16:63--76, 2007.

\bibitem{Morimoto}
M.~Morimoto.
\newblock {\em Analytic functionals on the sphere}.
\newblock American Mathematical Society, 1962.

\bibitem{Nussbaum}
A.~E. Nussbaum.
\newblock Quasi-analytic vectors.
\newblock {\em Archiv der Mathematik}, 6:179--191, 1965.

\bibitem{Pakes}
A.~G. Pakes, W.~Hung, and J.~Wu.
\newblock Criteria for the unique determination of probability distributions by
  moments.
\newblock {\em Australian and New Zealand Journal of Statistics}, 43:101--111,
  2001.

\bibitem{Pedersen}
H.~L. Pedersen.
\newblock On {K}rein's theorem for indeterminacy of the classical moment
  problem.
\newblock {\em J. Approx. Theory}, 95:90--100, 1998.

\bibitem{Petersen}
L.~C. Petersen.
\newblock On the relation between the multidimensional moment problem and the
  one-dimensional moment problem.
\newblock {\em Mathematica Scandinavia}, 51:361--366, 1982.

\bibitem{Rudin2}
W.~Rudin.
\newblock {\em Real and complex analysis}.
\newblock McGraw-Hill, 1973.

\bibitem{RudinCn}
W.~Rudin.
\newblock {\em Function theory in the unit ball of $\C^n$}.
\newblock Springer, 1980.

\bibitem{SZ}
H.~Salzmann and K.~Zeller.
\newblock Singularit\"aten unendlich oft differenzierbarer funktionen.
\newblock {\em Mathematische Zeitschrift}, 62:354--367, 1955.

\bibitem{SV}
J.~Schmets and M.~Valdivia.
\newblock On the extent of (non) quasi-analytic classes.
\newblock {\em Archiv der Mathematik}, 56:593--600, 1991.

\bibitem{shabat}
B.~V. Shabat.
\newblock {\em Introduction to complex analysis: functions of several
  variables, Part 2}.
\newblock AMS, 1992.

\bibitem{ST}
J.~Shohat and J.~Tamarkin.
\newblock {\em The problem of moments}.
\newblock American Mathematical Society, 1970.

\bibitem{JS}
J.~Stoyanov.
\newblock {\em Counterexamples in probability}.
\newblock John Wiley \& Sons, 1997.

\bibitem{Stoyanov}
J.~Stoyanov.
\newblock Krein condition in probabilistic moment problems.
\newblock {\em Bernoulli}, 6:939--949, 2000.

\bibitem{StoyanovLin}
J.~Stoyanov and G.~D. Lin.
\newblock Hardy's condition in the moment problem for probability
  distributions.
\newblock {\em Theory Probab. Appl.}, 57:811--820, 2012.

\bibitem{Tauvel}
P.~Tauvel.
\newblock {\em Analyse complexe pour la licence 3}.
\newblock Dunod, 2006.

\bibitem{Train}
K.~Train.
\newblock {\em Discrete choice methods with simulation}.
\newblock Cambridge, 2009.

\bibitem{Vytlacil}
E.~Vytlacil.
\newblock Independence, monotonicity, and latent index models: An equivalence
  result.
\newblock {\em Econometrica}, 70:331--341, 2002.

\bibitem{wu2010varying}
Y.~Wu, J.~Fan, H.-G. M{\"u}ller, et~al.
\newblock Varying-coefficient functional linear regression.
\newblock {\em Bernoulli}, 16(3):730--758, 2010.

\bibitem{Young}
R.~Young.
\newblock {\em An introduction to nonharmonic {F}ourier series}.
\newblock Academic press, 2001.

\end{thebibliography}

\end{document}